\documentstyle[12 pt]{article}
\def\endpf{\hbox{\vrule height1.5ex width.5em}}
\begin{document}
\title{ Identification of Parameters through the  Approximate
Periodic Solutions of a Parabolic System}
\author{  Ling Lei \\  Department of Mathematics, Zhejiang University,\\
 Hangzhou, 310027, P.R.China}
\date{}
\maketitle
{\bf Abstract.}    This work is concerned with the identification
problem for what we call the perturbation term or error term in a
parabolic partial differential equation, through its approximate
periodic solutions. The observation is made over a subregion of
the physical domain. The existence and uniqueness problem of the
approximate periodic solutions is  studied in the first part of
the paper.   A solution to the identification problem is  given in
the second part of the paper.  The main ingredients to be used
include  the classical Garlerkin method and the more recently
developed Carleman estimates for a parabolic system.
\\

{\bf Key words.}
 Identification of parameter, Garlerkin method, Carleman
inequality, approximate periodic solution, parabolic equation.
\\

{\bf AMS subject classification.     }  35K99, 93A99.\\ \vskip 1cm

\section{Introduction}
\hspace*{0.5 cm}Let $\Omega\subset {\bf R}^n$ be a bounded domain
with $C^2$-smooth boundary $\partial\Omega$ and let
$\omega\subset\Omega$ be a subdomain. Write $Q=\Omega \times
(0,T)$ with $T>0$ and  write $Q^{\omega}= \omega \times (0,T) $.
Consider the following parabolic equation:
$$
 \left\{\begin{array}{ll}
\displaystyle{\frac{\partial u}{\partial t}(x,t)}+Lu(x,t)=f(x,t),
\;\;\;\;\;\;\;\;\;\;\;\;\;\;\;\;\;\;\;\;\;&
\mbox{in }\;\;Q=\Omega \times (0,T),\\
 u(x,t)=0, & \mbox{on }\;\; \Sigma=\partial \Omega \times (0,T),\\
\end{array}\right.
\eqno{(1.1)}
$$
where $$Lu(x,t)=L_0 u(x,t)+ e(x,t)u(x,t),$$
$$
L_0u(x,t) = - \sum ^n _{i,j=1} D_j (a^{ij}(x) D_i u(x,t)) + \sum
^n_{i=1} b^i (x) D_i u(x,t) -\sum ^n_{i=1} D_i(b^i(x)u(x,t))
+c(x)u(x,t).
$$
Here and in what follows, we write $D_j=\frac{\partial}{\partial
x_j}$. We also  use the standard summation convection. Namely,
repeated indices imply summation from 1 to $n$. Throughout of the
paper, we make the following regularity assumptions for the
coefficients:

\noindent (I):  $a^{ij}(x) \in
Lip(\overline{\Omega}),\;a^{ij}(x)=a^{ji}(x),\; \mbox{and}\;
\lambda ^*|\xi|^2 \leq a^{ij}(x) \xi _i \xi_j \leq
\displaystyle{\frac{1}{\lambda^*}} |\xi|^2 ,\; \mbox{for}\; \xi
\in {\bf R}^n $ with $\lambda^*$  a  certain positive constant;

\noindent (II):  $b^i(x) \in Lip(\overline{\Omega}),\; c(x) \in
L^\infty(\Omega),\;\mbox{and}\;
 e(x,t) \in L^\infty (0,T;L^q(\Omega))\ \hbox{with }\ q
>\frac{n+2}{2}$ and $f(x,t)\in L^2(Q)$.\\

In many applications, one often encounters various problems, such
as the inverse problem and the Pontryagin maximum principle,
related to the periodic solutions of (1.1). (See [11] [15] [16]
[17], etc).  Here, we recall that a periodic solution of (1.1) is
a solution satisfying  the following condition:
$$
\begin{array}{ll}
u(x,0)=u(x,T) & \mbox{in }\;\;\Omega.
\end{array}
\eqno{(1.2)}
$$
In  (1.1),  the coefficients of the principal part $L_0$ of the
operator $L$ is $t$-independent. However, the one for what we call
the perturbation term or the error term $e(x,t)$ may well depend
on the time variable $t$. It is known that  for (1.1), when the
operator $L$ is not positive (for instance, when $e(x,t)$ takes
negative values), the periodic solution of (1.1)  may not exist
for a generic choice of $f(x,t)$. (See Example 3.4 in Section 3).
Namely, adding an error term with coefficient $e(x,t)$ to the
system may well destroy the periodicity  of certain solutions even
if $L_0$ is a positive operator. However, as we will show, the
system always possesses solutions with certain approximate
periodicity. This makes it a natural   problem to consider the
inverse problem, the Pontryagin problem, and many others, for
(1.1) through a certain family of solutions with  approximate
periodicity. In this paper, we will make an effort towards this
study by introducing the concept of  approximate periodic
solutions through the principle part $L_0$ of $L$. We then  study
the existence and uniqueness of such solutions and use them to
identify the error term through the observation of solutions over
$\omega$.

We next introduce the concept of $\mathcal{K}$-approximate
periodic solutions of (1.1), where $\mathcal{K}$ is a non-negative
integer.

First, we notice that $L_0$ is a symmetric operator.  Consider the
eigenvalue problem of $L_0$:
$$
\left\{
\begin{array}{ll}
L_0 v(x)=\lambda v (x),\\
v(x)|_{\partial \Omega} =0.
\end{array}\right.
\eqno{(1.3)}
$$
It is well-known (see [8]) that (1.3) has  a complete set of
eigenvalues $\{ \lambda _j\}^\infty _{j=1}$ with the associate
eigenvectors $\{X_j(x)\}^\infty _{j=1}$ such that $L_0
X_j(x)=\lambda _j X_j (x),\;-\infty<\lambda_1\leq\lambda_2
\leq\cdots\leq\lambda_j\leq\cdots<\infty, \lim_{j\rightarrow
\infty}\lambda_j=\infty, X_j(x)\in H^1_0(\Omega)$.
 Choose $\{X_j(x)\}^\infty _{j=1}$
so that it serves an orthonormal  basis of $L^2(\Omega)$. Therefore,
$\forall u(x,t)\in L^2(Q)$, we have
$u(x,t)=\sum^{\infty}_{j=1}u_j(t)X_j(x)$, where
$u_j(t)=\displaystyle{\int_{\Omega}}u(x,t)X_j(x)dx\in L^2(0,T)$.

\bigskip
{\bf Definition 1.1.} {\it We call $u(x,t)$ is a
$\mathcal{K}$-approximate periodic solution of (1.1) with respect
to its principal part $L_0$ if \\
(a):  $u \in C([0,T];L^2(\Omega))\cap L^2(0,T;H^1_0 (\Omega ))$ is
a weak solution of (1.1);\\
(b): $ u\in \mathcal{S}_{ \mathcal{K}}$ $ = \{u\in
C([0,T];L^2(\Omega));u_j(0)=u_j(T)\;for\;j \ge {\mathcal K}+1$ , $
u_j(t)=\displaystyle{\int_\Omega } u(x,t)X_j (x)dx \}.$}

\bigskip
Here, we recall that $u(x,t)$ is said to be a weak solution of (1.1)
with the initial value $u(x,0)=\psi(x)$ if
$u\in  C([0,T];L^2(\Omega))\cap L^2(0,T;H^1_0 (\Omega ))$ and
 for any testing function
$\varphi \in H^{1,1}_0(Q)=\{ h \in L^2 (Q);
\partial _t h \in L^2(Q),D_i h \in L^2(Q)\;for\; all\; i=1,2 \cdots n,\;
h(x,t)|_{\partial\Omega} =0 \}$, we have, for all $t$,
$$
(u(\cdot,t),\varphi (\cdot,t))
 -\displaystyle{\int ^t_0}(u,\varphi _\tau)d\tau
 +\displaystyle{\int ^t_0 }(Lu,\varphi) d\tau
=(\psi,\varphi(\cdot,0)) +\displaystyle{\int ^t_0}(f,\varphi
)d\tau.$$ When $\mathcal{K}$ $ = 0$, we will always  regard $\sum
^{0} _{j=1} = 0$.  Hence,  a 0-approximate periodic solution of
(1.1) is a regular periodic solution.

It should be mentioned that in the above
definition, we need only to assume that $u(x,t)\in L^\infty
(0,T;L^2(\Omega )) \cap L^2(0,T;H^1_0(\Omega ))$ to start with. Then it
holds automatically that $u\in C([0,T];L^2(\Omega))$ (see Chapter 3 of [5]).
Also, we notice that

$ u\in \mathcal{S}_{\cal K} \Longleftrightarrow$
$(u-\sum^{\cal{K}} _{j=1} (u,X_j)X_j)(x,0)
=(u-\sum^{\cal{K}} _{j=1} (u,X_j)X_j)(x,T).$

In the above formula and in what follows, we write
$(u(\cdot,t),\varphi(\cdot ,t))=\displaystyle{\int_\Omega}
u(x,t)\varphi(x,t)dx$, $(u(\cdot,t),u(\cdot ,t))=\|u(\cdot,t)\|^2$,
and we denote $u_t$  for the derivative of $u(x,t)$ with respect to $t$.

Our first result of  this paper can be stated  as follows:

\bigskip
{\bf Theorem 1.2}. {\it Let $e(x,t)\in$ ${\mathcal{M}}_q$ $=
\{e(x,t) \in L^{\infty}(0,T;L^q (\Omega));$ $\mbox{ess sup}_{t\in
(0,T)} \|e(x,t)\|_{L^q (\Omega)}$ $ \leq M,q> \frac{n+2}{2},M\;
\mbox{is a constant}\}$. Then, there exists an integer ${\mathcal
K}_0\equiv{\mathcal K}_0$ $(L_0,M,\Omega)$ $\geq 0$ such that for
any ${\mathcal K}\geq{\mathcal K}_0$ and any initial value
$a_I=(a_1,a_2,\cdots,a_{\mathcal{K}})\in {\bf R}^{\mathcal K}$, we
have a unique solution to the following equation:
$$
 \left\{\begin{array}{ll}
\displaystyle{\frac{\partial u(x,t)}{\partial t}}+L_0u(x,t)+e(x,t)u(x,t)=f(x,t),
\;\;\;\;\;\;\;\;\;\;\;\;\;\;\;\;\;\;\;\;\;&
\mbox{in }\;\;Q,\\
 u(x,t)=0, & \mbox{on }\;\; \Sigma,\\
 (u(x,0),X_j(x))=a_j, & \mbox{for }\;j\leq \mathcal{K},\\
u\in \mathcal{S}_{ \mathcal{K}}.
\end{array}\right.
\eqno{(1.4)}
$$
Moreover, for such a solution $u(x,t)$, we have the following
energy estimate:
$$
\begin{array}{ll}
\sup_{t\in[0,T]}\|u(\cdot,t)\|^2 +\displaystyle{\int^T_0}\|\nabla
u(\cdot,t)\|^2 dt \leq C(L_0,M,\Omega) (|a_I|^2
+\displaystyle{\int_Q}f^2dxdt).
\end{array}
\eqno (1.5) $$}
\bigskip

The second part of this work is to study an inverse problem. We
will identify $(e(x,t),a_I)$ from ${\mathcal{M}}_q\times {\bf
R}^{\mathcal{K}}$ via the observation of solutions for (1.4) on
the subdomain $\omega \subset \Omega$. More precisely, we shall
study the following identification problem:
\bigskip

$Problem \ (P)\ \ \ \hbox{Find   the  minimum  value  of}\;
\displaystyle{\int_{Q^{\omega}}}|u-\widetilde{u}|^{2}dxdt \mbox{
for }(e(x,t),a_I)\in {\mathcal{M}}_q \times {\bf R}^{\mathcal{K}}$
with $u $ satisfying
$$
 \left\{\begin{array}{ll}
\displaystyle{\frac{\partial u(x,t)}{\partial t}}+L_0u(x,t)+e(x,t)u(x,t)=f(x,t),
\;\;\;\;\;\;\;\;\;\;\;\;\;\;\;\;\;\;\;\;\;&
\mbox{in }\;\;Q,\\
 u(x,t)=0, & \mbox{on }\;\; \Sigma,\\
 (u(x,0),X_j(x))=a_j, & \mbox{for }\;j\leq \mathcal{K},\\
u\in \mathcal{S}_{ \mathcal{K}},
\end{array}\right.
$$
where $ \widetilde{u} \in L^2(Q^\omega)$ is a
given function.
\bigskip

Making use of Theorem 1.2 and  the Carleman inequality established
in [2] [11]  [16], etc, we are able to prove  the existence of
solutions to problem $(P)$. Our second main result  can be stated
as follows:

\bigskip
{\bf Theorem 1.3}. {\it  Let $\mathcal K$ be as in Theorem 1.2. Then
there exist an  $e^*(x,t)\in
{\mathcal{M}}_q$ and $a^*_I \in {\bf R}^{\mathcal{K}}$ such that
$$
\displaystyle{\int_{Q^{\omega}}}|u(e^*,a^*_I;x,t)-
\widetilde{u}|^{2}dxdt= \inf_{(e,a_I)\in {\mathcal{M}}_q\times
{\bf R}^{\mathcal{K}}}
\displaystyle{\int_{Q^{\omega}}}|u(e,a_I;x,t)-
\widetilde{u}|^{2}dxdt.
$$
Here $ \widetilde{u} \in L^2(Q^\omega)$ is a given function and
$u(e,a_I;x,t)$ is the solution of equation (1.4) with error
coefficient $e(x,t)$ and $(u(e,a_I;\cdot,0),X_j(\cdot))=a_j$ for
$j\le {\mathcal K}$, where $a_I=(a_1,\cdots,a_{\mathcal K})\in
{\bf R}^{\mathcal{K}}$.}
\medskip

Theorem 1.3 can be immediately used to give the following slightly
more general result:

\bigskip
{\bf Corollary 1.4}. {\it  Let $k$ be a non-negative integer. Then
there exist an  $e^*(x,t)\in {\mathcal{M}}_q$ and $a^*_I \in {\bf
R}^{{k}}$ such that
$$
\displaystyle{\int_{Q^{\omega}}}|u^*-
\widetilde{u}|^{2}dxdt= \inf_{(e,a_I)\in {\mathcal{M}}_q\times
{\bf R}^{k}, u\in U(e,a_I;x,t)}
\displaystyle{\int_{Q^{\omega}}}|u-
\widetilde{u}|^{2}dxdt.
$$
Here $ \widetilde{u} \in L^2(Q^\omega)$ is a given function and
$U(e,a_I;x,t)$ is the set of solutions of the following equation,
which we assume to be non-empty:
 $$
 \left\{\begin{array}{ll}
\displaystyle{\frac{\partial u(x,t)}{\partial t}}+L_0u(x,t)+e(x,t)u(x,t)=f(x,t),
\;\;\;\;\;\;\;\;\;\;\;\;\;\;\;\;\;\;\;\;\;&
\mbox{in }\;\;Q,\\
 u(x,t)=0, & \mbox{on }\;\; \Sigma,\\
 (u(x,0),X_j(x))=a_j,\ a_I=(a_1,\cdots,a_k),\  & \mbox{for }\;j\leq k,\\
u\in {\mathcal S}_{k}.
\end{array}\right.
\eqno (1.6)$$}

Similarly, $u^*$ satisfies the following equation:
$$
 \left\{\begin{array}{ll}
\displaystyle{\frac{\partial u^*(x,t)}{\partial t}}+L_0u^*(x,t)+e^*(x,t)u^*(x,t)=f(x,t),
\;\;\;\;\;\;\;\;\;\;\;\;\;\;\;\;\;\;\;\;\;&
\mbox{in }\;\;Q,\\
 u^*(x,t)=0, & \mbox{on }\;\; \Sigma,\\
 (u^*(x,0),X_j(x))=a^*_j,\ a^*_I=(a^*_1,\cdots,a^*_k),\  & \mbox{for }\;j\leq k,\\
u^*\in {\mathcal S}_{k}.
\end{array}\right.
$$
\bigskip

Notice that in Corollary 1.4, (1.6) may have a family of different
solutions.

\medskip
It is not clear to us if  the uniqueness property for $e^*(x,t)$
in Theorem 1.3 holds. However, if one fixes $e(x,t)$ and tries to
identify $a_I$ through problem $(P)$, then the uniqueness of $a_I$
is indeed guaranteed as the following theorem shows:

\bigskip
{ \bf Theorem 1.5}. {\it  Under the same notation as in Theorem 1.3,
there exists a unique $a^*_I \in
{\bf R}^{\mathcal{K}}$ such that
$$
\displaystyle{\int_{Q^{\omega}}}|u(e,a^*_I;x,t)-
\widetilde{u}|^{2}dxdt= \inf_{a_I\in {\bf R}^{\mathcal{K}}}
\displaystyle{\int_{Q^{\omega}}}|u(e,a_I;x,t)-
\widetilde{u}|^{2}dxdt.$$}

\medskip
System (1.1) models a large class of physical processes, where
$u(x,t)$ represents the temperature or other physical quantity.
The identification problems associated with system (1.1) with
initial condition $u(x,0)=u_0(x)$, where $u_0(x)$ is a given
function, were studied by many authors. See [1] [4] [6] [7] and
[13], where the observations are taken in the whole domain
$\Omega$. However, in many applications, one may only be able to
measure the quantity  on a subdomain $\omega\subset\Omega$ and
does not have enough information about the initial value. One may
still be asked to determine the error influence $e(x,t)$ in the
physical process through  the approximate value of the  solutions
over $\omega$. When the approximate value comes from  approximate
periodic solutions, then our results of the present paper can be
directly applied. Notice that  this is an inverse problem. For the
direct problem, one is asked to determine the value of the
solution to (1.1) for a given
$e(x,t)$ and $a_I$.\\
\hspace*{0.5 cm}Since our observation is taken in a subdomain
$\omega\subset\Omega$, we can not apply the method employed in the
work mentioned above to answer $(P)$. Our key ingredients in this
paper to get the existence of  solutions for $(P)$  are  the
energy estimate (1.5) and  the Carleman inequality. There have
been many papers written on the related subjects in recent years.
Here we would like to mention [1-7] [14] [16] [17],  and the
reference therein, to name a few.\\
\hspace*{0.5 cm}The paper is organized as follows. In Section 2,
we prove the existence and uniqueness of the solution to system
(1.4). In Section 3, we obtain the existence of the identification
problem $(P)$ by proving Theorem 1.3.

\medskip

 This work is a continuation of  [14], where  very special cases of
 the results in this paper were studied. It should be mentioned
that the paper is largely motivated by two  papers of G. Wang and
L. Wang [16] [17].

\bigskip

\section{The existence and uniqueness of the solution}
\hspace*{0.5 cm}In this section, we prove the existence and
uniqueness of the solution to system (1.4). We will use the
 Galerkin method for constructing  solutions in the
$\mathcal{S}_{ \mathcal{K}}$-space for  the following equation
introduced  in Section 1:
$$
\left\{
\begin{array}{ll}
\displaystyle{\frac{\partial u(x,t)}{\partial t}}+L_0u(x,t)+e(x,t)u(x,t)=f(x,t),
\;\;\;\;\;\;\;\;\;\;\;\;\;\;\;\;\;\;\;\;\;&
\mbox{in }\;\;Q,\\
 u(x,t)=0, & \mbox{on }\;\; \Sigma.\\
\end{array}\right.
\eqno{(2.1)}
$$
Recall that $L_0 X_j(x)=\lambda _j X_j(x),\lambda
_j\rightarrow\infty$, $X_j(x)\in H^1_0(\Omega)$. Let
$G_N=\{g(x,t)\in L^2(Q);g(x,t)=\sum^N_{j=1} g_j (t)X_j(x) ,$
$g_j(t)\in L^2(0,T)\}$. We first look for an approximate solution
$u^N(x,t)$ of (2.1) in the $G_N$-space, which also has the
$\mathcal{K}$-approximate periodicity   as defined before. Here
$\mathcal{K}$ depends only on the $L_0,\ M, \Omega$ and will be
determined   later. $N$ is always assumed to be sufficiently large
($N>>{\mathcal{K}}$).\\
\hspace*{0.5 cm}Write $L=L_0+e$. Assume $u^N = \sum^N_{j=1}
u^N_j(t)X_j(x)$  such that $\partial_t u^N = -Lu^N +f$ has 0
projection to $G_N$ in the following sense:
$$
(\partial_t u^N +Lu^N -f,\;\varphi)=0,\;\;\;\;\;
\mbox{for}\;0<t<T\;\mbox{and any}\;\varphi \in G_N.
$$
Letting $\varphi = X_j$ for $j=1,2,\cdots ,N$, we get the
following system of ordinary differential equations:
$$\begin{array}{ll}
\displaystyle{\frac{du^N_j (t)}{dt}}+\sum B_{kj}(t)u^N_k
(t)=f_j(t),\; j=1,2,\cdots,N,\\
\mbox{where}\; B_{kj}(t)=(LX_k,X_j)=\displaystyle{\int_\Omega}
LX_k \cdot X_j
dx,\;f_j(t)=(f,X_j)=\displaystyle{\int_\Omega}f(x,t)X_j(x)dx.
\end{array}
$$
We put the following condition on $u_j^N$:
$$
u^N_j(0)=u^N_j(T)\;\;\mbox{for}\; j >
\mathcal{K}\;\mbox{with}\;\mathcal{K} \;\mbox{ independent of $N$
and being determined later}.
$$
Consider the following system of ordinary differential equations:
$$\left\{\begin{array}{ll}
\displaystyle{\frac{du^N_j (t)}{dt}}+\sum B_{kj}(t)u^N_k
(t)=0,\ j=1,\cdots, N,\\
u^N_I(0)=0,\\
u^N_{II}(0)=a^N_{II}\in {\bf R}^{N-\mathcal{K}}.
\end{array}\right.
\eqno{(2.2)}
$$
Here and in what follows,
$$\begin{array}{ll}
u^N_I(t)=(u^N_1(t),u^N_2(t),\cdots , u^N_{\mathcal{K}}(t)),\\
u^N_{II}(t)=(u^N_{{\mathcal{K}}+1}(t),u^N_{{\mathcal{K}}
+2}(t),\cdots , u^N_N(t)).
\end{array}
$$
\hspace*{0.5 cm}{\bf Lemma 2.1.} {\it Let
$u^N(x,t)=\sum^{N}_{j=1}u^N_j(t)X_j(x)$ be the solution of (2.2).
There exists an integer $\mathcal{K}$ depending only on $L_0,
M,\Omega$  such that for any fixed $N>\mathcal{K}$, the operator:
$$
J:\;{\bf R}^{N-\mathcal{K}}\longmapsto
{\bf R}^{N-\mathcal{K}},\;J(a^N_{II})=u^N_{II}(T),
$$ is contractive. Namely,
$$
|J(a^N_{II})| \leq \mu |a^N_{II}|\;\;\mbox{with}\; \mu\;
\mbox{fixed and }0\leq \mu <1.
$$
Here and in what follows, we always assume that $\mbox{ess sup}_{t\in (0,T)}
\|e(x,t)\|_{L^q (\Omega)}\leq M$.}\\

\medskip
For the proof of Lemma 2.1, we need the following
claim:\\\\
\hspace*{0.5 cm}{\bf Claim 2.2.}\ {\it For any $v_1,v_2 \in
H^1_0(\Omega)$, there are constants $C(\Omega,q)$ depending only
on $\Omega,q$ and $C_s(\varepsilon)$ and $C_l(\varepsilon)$,
depending only on $\varepsilon$ with $C_s(\varepsilon)\rightarrow
0$ as $\varepsilon\rightarrow 0$, $C_l(\varepsilon)\rightarrow
\infty $ as $\varepsilon\rightarrow 0$, such that
$$\begin{array}{ll}
\displaystyle{\int_\Omega} |e(x,t)v_1(x)v_2(x)|dx &\leq
C(\Omega,q)M \{C_l(\varepsilon)(\|v_1\|^2_{L^2(\Omega)} +
\|v_2\|^2_{L^2(\Omega)})\\
&+C_s(\varepsilon)(\|\nabla v_1\|^2_{L^2(\Omega)} +\|\nabla
v_2\|^2_{L^2(\Omega)})\}.
\end{array}
$$}

\medskip
{\it Proof of Claim 2.2.}\ By the Schwartz inequality, we need
only to prove the  claim in the case of $v_1 =v_2=v$.
$$
\begin{array}{ll}
\displaystyle{\int_\Omega} |e(x,t)|v^2dx &\leq \|e(\cdot,
t)\|_{L^q(\Omega)}\|v^2\|_{L^{q'}(\Omega)} \\
&\leq M \|v\|^2_{L^{2q'}(\Omega)},
\end{array}
$$
where $ q'=\frac{q}{q-1}$. Let $\alpha =\frac{1}{2q'} (n+2-nq')$.
Then
$\frac{1}{2q'}=\frac{\alpha}{2}+\frac{1-\alpha}{\frac{2(n+2)}{n}}$.
Next, by the H$\ddot{o}$lder inequality, we have
$$
I=\displaystyle{\int_\Omega} |e(x,t)|v^2dx \leq M
\|v\|^{2\alpha}_{L^2(\Omega)}\|v\|^{2(1-\alpha)}
_{L^{\frac{2(n+2)}{n}}(\Omega)}.
$$
By the Sobolev inequality,
$$
\|v\|_{L^{\frac{2(n+2)}{n}}(\Omega)}\leq C(\Omega,q)\|v\|_{H^1_0
(\Omega)}.
$$
Hence, $I\leq M C(\Omega,q)
(\|v\|^{\alpha}_{L^2(\Omega)}\|v\|^{1-\alpha}_{H^1_0 (\Omega)}
)^2$.\\

Now, by the following H$\ddot{o}$lder inequality:
 $$a\cdot b \leq
\varepsilon a^p + \frac{1}{\varepsilon^{\frac{p}{q}}} b^q$$
 with
$\frac{1}{p}+\frac{1}{q} =1$, we have
$$
\begin{array}{ll}
I &\leq M C(\Omega,q) (C_s(\varepsilon)\|v\|_{H^1_0
(\Omega)}+C_l(\varepsilon)\|v\|_{L^2(\Omega)} )^2 \\
& \leq M C(\Omega,q)(C_s(\varepsilon)\|v\|^2_{H^1_0
(\Omega)}+C_l(\varepsilon)\|v\|^2_{L^2(\Omega)} )\\
&\leq M C(\Omega,q)(C_s(\varepsilon)\|\nabla
v\|^2_{L^2(\Omega)}+C_l(\varepsilon)\|v\|^2_{L^2(\Omega)} ).
\end{array}
$$
Here and in what follows, $C_s(\varepsilon),C_l(\varepsilon)$
stand for small and large constant depending only on
$\varepsilon$, which may be different in different contexts. The
proof of the claim is complete. $\endpf$

\medskip
By Claim 2.2,  we have
$$
\begin{array}{ll}
|B_{kj}(t)|&=|(L_0 X_k, X_j)+(e(x,t)X_k, X_j)|\\
&=|\lambda _j \delta ^k_j +(e(x,t)X_k, X_j)|\\
&\leq |\lambda _j \delta ^k_j| +|C_l(\varepsilon)+
C_s(\varepsilon)(\|\nabla X_k\|^2_{L^2(\Omega)} +
\|\nabla X_j\|^2_{L^2(\Omega)})|,\\
|f_j(t)| &\leq (\displaystyle{\int _{\Omega}}f^2
(x,t)dx)^{\frac{1}{2}}( \int _{\Omega}X^2_j
dx)^{\frac{1}{2}}\leq\|f\|_{L^2(\Omega)}.
\end{array}
$$
In particular, we conclude that any solution of the initial value
problem of

$$
\left\{\begin{array}{ll} \displaystyle{\frac{du^N_j (t)}{dt}}+\sum
B_{kj}(t)u^N_k
(t)=f_j(t),\ j=1,\cdots,N,\\
(u^N_1(0),u^N_2(0),\cdots,u^N_N(0))=a^N \in {\bf R}^N,
\end{array}\right.
\eqno{(2.3)}
$$
is absolutely continuous over $[0,T]$.

\bigskip
{\it Proof of Lemma 2.1.}\ Multiplying $2u^N_j(t)$ to the first
equation of (2.2) and summing up with respect to j from 1 to $N$,
we get
$$
\frac{d(\|u^N(\cdot,t)\|^2)}{dt}+2(L_0
u^N,u^N)+2(e(x,t)u^N,u^N)=0.$$ As before, we use $\| \cdot\|$ to
denote  the usual $L^2(\Omega)$-norm. After some calculation
involving the Green formula, we have the following G\"arding
inequality (see [8]):
$$
\frac{2}{\lambda^*}\|\nabla u^N\|^2  +C\|u^N\|^2\geq(L_0u^N,u^N)
 \geq \frac{\lambda^*}{2}\|\nabla u^N\|^2  - C\|u^N\|^2.
 \eqno{(2.4)}
$$
By (2.4) and Claim 2.2, we obtain
$$
\frac{d(\|u^N(\cdot,t)\|^2)}{dt}+\lambda^*\|\nabla u^N\|^2
-C\|u^N\|^2-C_s(\varepsilon)\|\nabla u^N\|^2-C_l(\varepsilon)
\|u^N\|^2 \leq 0.
$$
We choose $\varepsilon$ such that
$C_s(\varepsilon)<\frac{\lambda^*}{2}$. Then we get, for a large
constant $C_l$,
$$
\frac{d(\|u^N(\cdot,t)\|^2)}{dt} +\frac{\lambda^*}{2}\|\nabla
u^N\|^2 -C_l(\varepsilon) \|u^N\|^2\leq 0.
$$
Applying the Gronwall inequality, we have
$$
\frac{d}{dt}(\|u^N(\cdot,t)\|^2 e^{-C_l(\varepsilon) t}) +
\frac{\lambda^*}{2} e^{-C_l(\varepsilon) t}\|\nabla u^N\|^2 \leq
0,
$$
$$
\|u^N(\cdot,t)\|^2 + \displaystyle{\int^t_0}\|\nabla
u^N(\cdot,\tau)\|^2d \tau  \leq C
\|u^N(\cdot,0)\|^2,\;\;\;\;\;\forall t \in [0,T].
$$
In particular,
$$
\|u^N(\cdot,T)\|^2\leq C \|u^N(\cdot,0)\|^2 =C|a^N_{II}|^2,\;
\mbox{and}\;\displaystyle{\int^T_0}\|\nabla u^N(\cdot,t)\|^2dt
\leq C|a^N_{II}|^2. \eqno{(2.5)}
$$
Next, multiplying $ 2u^N_j(t)$ to the first equation of (2.2) and
summing up with respect to $j$ from $\mathcal{K}+ \mbox{1}$ to
$N$, and letting
$$
u^{N,II}(x,t)=\sum ^N_{j={\mathcal{K}}+1} u^{N}_{j}(t)X_j(x),
$$
then we get
$$
\frac{d(\|u^{N,II}(\cdot,t)\|^2)}{dt}+2(Lu^N,u^{N,II})=0.
$$
Notice that
$$
\begin{array}{ll}
(Lu^N,u^{N,II})&=(L_0 u^N,u^{N,II})+(e(x,t)u^N,u^{N,II}) \\
&=(L_0 u^{N,II},u^{N,II})+(e(x,t)u^N,u^{N,II})\\
&=\sum^N_{j={\mathcal{K}}+ 1} \lambda_j (u^N_j(t))^2
+(e(x,t)u^N,u^{N,II}).\\
\end{array}
$$
By Claim 2.2 and (2.4), we have
$$
\begin{array}{ll}
|(e(x,t)u^N,&u^{N,II})|\leq C_l(\varepsilon)(\|u^N\|^2 +
\|u^{N,II}\|^2)  + C_s(\varepsilon)(\|\nabla u^N\|^2 + \|\nabla
u^{N,II}\|^2)\\
&\leq C_l(\varepsilon)\|u^N\|^2 +C_s(\varepsilon)\|\nabla
u^N\|^2+C_s(\varepsilon)\{\frac{2C}{\lambda^*}\|u^{N,II}\|^2+\frac{2}{\lambda^*}(L_0
u^{N,II},u^{N,II})\}\\
&\leq C_l(\varepsilon)\|u^N\|^2 +C_s(\varepsilon)\|\nabla
u^N\|^2+C_s(\varepsilon)\{\frac{2C}{\lambda^*}\|u^{N}\|^2+\frac{2}{\lambda^*}
\sum^N_{j={\mathcal{K}}+1} \lambda_j (u^N_j(t))^2\}\\
&\leq C_l(\varepsilon)\|u^N\|^2 +C_s(\varepsilon)\|\nabla
u^N\|^2+C_s(\varepsilon)\{\frac{2C}{\lambda^*}\|u^{N}\|^2+\frac{2}{\lambda^*}
\sum^N_{j=1} \lambda_j (u^N_j(t))^2\}\\
&\leq C_l(\varepsilon)\|u^N\|^2 +C_s(\varepsilon)\|\nabla
u^N\|^2+C_s(\varepsilon)\{\frac{2C}{\lambda^*}\|u^{N}\|^2+\frac{2}{\lambda^*}
(L_0u^N,u^N)\}\\
&\leq C_l(\varepsilon)\|u^N\|^2 +C_s(\varepsilon)\|\nabla
u^N\|^2+C_s(\varepsilon)(\frac{4C}{\lambda^*}\|u^{N}\|^2+
\frac{4C}{\lambda^{*2}}\|\nabla u^{N}\|^2)\\
&\leq C_l(\varepsilon)\|u^N\|^2 +C_s(\varepsilon)\|\nabla u^N\|^2.
\end{array}
$$
Notice that to get the last inequality, we applied the other part
of the G\"{a}rding estimate. We have
$$
\frac{d(\|u^{N,II}(\cdot,t)\|^2)}{dt} + 2\lambda
_{\mathcal{K}}\|u^{N,II}\|^2 -C_l(\varepsilon)\|u^N\|^2
-C_s(\varepsilon)\|\nabla u^N\|^2 \leq 0.
$$
By the Gronwall inequality, we get
$$
\begin{array}{ll}
e^{2\lambda _{\mathcal{K}}
t}\|u^{N,II}(\cdot,t)\|^2-\|u^{N,II}(\cdot,0)\|^2 &\leq
C_l(\varepsilon)\displaystyle{\int ^t_0}e^{2\lambda _{\mathcal{K}}
\tau }\|u^{N}\|^2 d\tau\\
 &+C_s(\varepsilon)\displaystyle{\int
^t_0}e^{2\lambda _{\mathcal{K}} \tau }\| \nabla u^{N}\|^2
d\tau,\;\;\forall t \in [0,T].
\end{array}
$$
We obtain
$$
\begin{array}{ll}
e^{2\lambda _{\mathcal{K}}
T}\|u^{N,II}(\cdot,T)\|^2-\|u^{N,II}(\cdot,0)\|^2 &\leq
C_l(\varepsilon)\displaystyle{\int ^T_0}e^{2\lambda _{\mathcal{K}}
t}\|u^{N}\|^2 dt\\
 &+C_s(\varepsilon)\displaystyle{\int
^T_0}e^{2\lambda _{\mathcal{K}} t}\| \nabla u^{N}\|^2 dt,
\end{array}
$$

$$
\begin{array}{ll}
\|u^{N,II}(\cdot,T)\|^2 & \leq e^{ -2 \lambda _{\mathcal{K}}
T}\|u^{N,II}(\cdot,0)\|^2  +C_l(\varepsilon)\displaystyle{\int
^T_0}e^{2\lambda _{\mathcal{K}}
(t-T)}\|u^{N}\|^2 dt\\
 &+C_s(\varepsilon)\displaystyle{\int
^T_0}e^{2\lambda _{\mathcal{K}} (t-T)}\| \nabla u^{N}\|^2 dt\\
&\leq e^{ -2 \lambda _{\mathcal{K}} T} |a^N_{II}|^2
+C_l(\varepsilon)\cdot C \cdot |a^N_{II}|^2(\frac{1}{2\lambda
_{\mathcal{K}}} - \frac{e^{ -2 \lambda _{\mathcal{K}} T}
}{2\lambda _{\mathcal{K}}})+C_s(\varepsilon)\cdot C \cdot
|a^N_{II}|^2 .
\end{array}
$$
Now, we first choose $\varepsilon$ sufficient small such that
$C_s(\varepsilon)\cdot C <\frac{1}{4}$. Then we fix such an
$\varepsilon$ and fix a $\mathcal{K}\gg$ $1$ such that $
e^{-2\lambda_{\mathcal{K}}T}<\frac{1}{4}$ and
$C_l(\varepsilon)\cdot C \cdot (\frac{1}{2\lambda_{\mathcal{K}}}
-\frac{e^{ -2\lambda_{\mathcal{K}} T}} {2\lambda _{\mathcal{K}}})
<\frac{1}{4}$. (Apparently, the choice of such a ${\mathcal K}$
depends only on the operator $L_0,\ M, \ \Omega$.) We then obtain
$$
|u^{N}_{II}(T)|^2=\|u^{N,II}(\cdot,T)\|^2 \leq
\frac{3}{4}|a^N_{II}|^2. \eqno{(2.6)}
$$
Since $J(a^N_{II})=u^{N}_{II}(T)$, we see the proof of Lemma 2.1.
$\endpf$\\

Next, we prove the following proposition:

\medskip
{\bf Proposition 2.3.} {\it Let $\mathcal{K}$ be choose as above.
Then for any $a^N_I \in {\bf R}^{\mathcal{K}}$, there exists a
unique solution $u^N$ of the following mixed boundary value
problem:
$$\left\{\begin{array}{ll}
\displaystyle{\frac{du^N_j (t)}{dt}}+\sum B_{kj}(t)u^N_k
(t)=f_j(t), \ j=1,\cdots, N, \\
u^N_I(0)=a^N_I\in {\bf R}^{\mathcal{K}},\\
u^N_{II}(0)=u^N_{II}(T).
\end{array}\right.
\eqno{(2.7)}
$$}
\hspace*{0.5 cm}{\it Proof of Proposition 2.3.}\ Let
$u^N_*(x,t)=\sum^{N}_{j=1}u^{N}_{*,j}(t)X_j(x)$ be the solution of
the following system:
$$\left\{\begin{array}{ll}
\displaystyle{\frac{du^N_j (t)}{dt}}+\sum B_{kj}(t)u^N_k
(t)=f_j(t),\ j=1,\cdots, N, \\
(u^N_1(0),u^N_2(0),\cdots,u^N_N(0))=0.\\
\end{array}\right.
\eqno{(2.8)}
$$
Let $u^N_{a^N}(x,t)=\sum^{N}_{j=1}u^{N}_{a^N,j}(t)X_j(x)$ be the
solution of the following system:
$$
\left\{\begin{array}{ll} \displaystyle{\frac{du^N_j (t)}{dt}}+\sum
B_{kj}(t)u^N_k
(t)=0,\ j=1,\cdots, N,\\
(u^N_1(0),u^N_2(0),\cdots,u^N_N(0))=a^N=(a^N_I,a^N_{II}).\\
\end{array}\right.
\eqno{(2.9)}
$$
For a fixed $a^N_I \in {\bf R}^{\mathcal{K}}$, we define
$$\widetilde{J}_{a^N_I} :
{\bf R}^{N-\mathcal{K}} \longrightarrow {\bf R}^{N-\mathcal{K}}\;\mbox{such
that } \widetilde{J}_{a^N_I} (a^N_{II})=u^N_{*II}(T)+u^N_{a^N
II}(T),
$$where
$$
\begin{array}{ll}
u^N_{*II}(T)&=(u^N_{*,{\mathcal{K}}+1}(T), u^N_{*,{\mathcal{K}}+2}(T),\cdots, u^N_{*,N}(T)),\\
u^N_{a^N II}(T)&=(u^N_{a^N, {\mathcal{K}}+1}(T), u^N_{a^N,
{\mathcal{K}}+2}(T),\cdots,u^N_{a^N, N}(T)).
\end{array}
$$
Namely, $u^N_{* II}$ is the last $N-{\mathcal K}$ components of $u^N_{*}$ and
$u^N_{a^N II}$ is the last $N-{\mathcal K}$ components of $u^N_{a^N}$.
We have
$$
\begin{array}{ll}
 |\widetilde{J}_{a^N_I} (a^{N,1}_{II})-\widetilde{J}_{a^N_I}
 (a^{N,2}_{II})|&=|u^N_{a^{N,1}II}(T) -u^N_{a^{N,2}II}(T)|\\
 &=|J(a^{N,1}_{II} -a^{N,2}_{II})|\\
 &\leq \frac{\sqrt{3}}{2}|a^{N,1}_{II} -a^{N,2}_{II}|.
\end{array}
$$
Hence, $\widetilde{J}_{a^N_I}$ is a contractive map and has a
unique fixed point $\widetilde{a^N_{II}}$. Namely,
$\widetilde{J}_{a^N_I}(\widetilde{a^N_{II}})=\widetilde{a^N_{II}}$.\\
Then (2.7) has a solution $u^N(x,t)$ with
$u^N_{II}(0)=\widetilde{a^N_{II}}=u^N_{II}(T)$. The uniqueness
also follows from the uniqueness of the fixed point of
$\widetilde{J}_{a^N_I}$. The proof of the proposition is
complete. $\endpf$

\medskip
Next, we estimate $\widetilde{a^N_{II}}$, the unique fixed point
of $\widetilde{J}_{a^N_I}$.

\medskip
{\bf Proposition 2.4.} {\it Let $\widetilde{a^N_{II}}$ be as
above. Then $|\widetilde{a^N_{II}}|^2 \leq C (|a^N_I|^2
+\displaystyle{\int_Q}f^2dxdt)$, where $C$ depends only on  $L_0,
M, \Omega$.}

\medskip
{\it Proof of Proposition 2.4.}\ We know
$\widetilde{a^N_{II}}=u^N_{II}(T)=u^N_{*,II}(T)+u^N_{a^N ,II}(T)$.
We first estimate $u^N_{*,II}(T)$. Multiplying the first equation
of (2.8) by $2u^N_{*,j}(t)$ and summing up with respect to $j$
from 1 to $N$, we have
$$
\begin{array}{ll}
\displaystyle{\frac{d(\|u^N_*(\cdot,t)\|^2)}{dt}} + 2(Lu^N_*,
u^N_*)&=\sum^N_{j=1}f_j(t)u^N_{*j}(t)\\
&\leq \|f\|^2 + \|u^N_*(\cdot,t)\|^2.
\end{array}
$$
By using the Gronwall inequality as before, we get
$$
\sup _{t\in[0,T]}\|u^N_*(\cdot,t)\|^2
+\displaystyle{\int^T_0}\|\nabla u^N_*(\cdot,t)\|^2 dt\leq C
\int_Q f^2 dxdt.
$$
In particular, $\|u^N_*(\cdot,T)\|^2\leq C \displaystyle{\int_Q}
f^2 dxdt$. Hence,
$$|u^N_{*,II}(T)|^2=\|u^N_{*II}(\cdot,T)\|^2\leq C
\displaystyle{\int_Q} f^2 dxdt. \eqno{(2.10)}
$$
Next, we let $u^{N,1} (x,t), u^{N,2}(x,t)$ be the solution of the following system
(2.11) and (2.12), respectively,
$$
\left\{
\begin{array}{ll}
\displaystyle{\frac{du^N_j (t)}{dt}}+\sum
B_{kj}(t)u^N_k(t)=0, \ j=1,\cdots, N, \\
(u^N_1(0),u^N_2(0),\cdots,u^N_N(0))=(a^N_I,0),
\end{array}\right.
\eqno{(2.11)}
$$

$$
\left\{\begin{array}{ll} \displaystyle{\frac{du^N_j (t)}{dt}}+\sum
B_{kj}(t)u^N_k(t)=0,\ j=1,\cdots, N, \\
(u^N_1(0),u^N_2(0),\cdots,u^N_N(0))=(0,\widetilde{a^N_{II}}).
\end{array}\right.
\eqno{(2.12)}
$$
From the proof of Proposition 2.3 and Lemma 2.1, we get
$$
\begin{array}{ll}
|u^N_{a^N II}(T)|^2&= |u^{N,1}_{II}(T) +u^{N,2}_{II}(T)|^2 \\
&\leq C_l(\varepsilon)|a^N_I|^2+ (1+C_s(\varepsilon))|J(\widetilde{a^N_{II}})|^2\\
&\leq C_l|a^N_I|^2+ (1+C_s)\frac{3}{4}|\widetilde{a^N_{II}}|^2.\\
\end{array}
\eqno{(2.13)}
$$
By (2.10) and (2.13), we have
$$|\widetilde{a^N_{II}}|^2 \leq C
(|a^N_I|^2 +\displaystyle{\int_Q}f^2dxdt).
$$
The proof of the proposition is complete. $\endpf$\\

\medskip
{\bf Remark 2.5.}\ From the proof of Lemma 2.1, we conclude that
the value $\mathcal{K}$ depends only on $L_0$, $\Omega$ and $M$.
Namely, given $L_0$ and $M$ with $\mbox{ess sup}_{t\in (0,T)}
\|e(x,t)\|_{L^q (\Omega)} \leq M$, there is a ${\mathcal
K}_0\equiv{\mathcal K}_0$ $(L_0,M,\Omega)$ such that for any
$N>\mathcal K \geq$ ${\mathcal K}_{0}$ and
$a^N_I=(a^N_1,a^N_2,\cdots, a^N_{\mathcal{K}})$ $\in {\bf
R}^{\mathcal{K}}$, the following mixed value problem has a unique
absolutely continuous solution
$u^N(x,t)=\sum^N_{j=1}u^N_j(t)X_j(x)$ over $[0,T]$:
$$
\left\{
\begin{array}{ll}
\displaystyle{\frac{du^N_j (t)}{dt}}+\sum B_{kj}(t)u^N_k
(t)=f_j(t),\ j=1,\cdots,N, \\
u^N_I(0)=a^N_I,\\
u^N_{II}(0)=u^N_{II}(T).
\end{array}\right. \eqno (2.14)
$$
Moreover, write $a^N_{II}=u^N_{II}(0)=u^N_{II}(T)$, we have the
estimate:
$$
\begin{array}{ll}
|a^N_{II}|^2 \leq C (|a^N_I|^2 +\displaystyle{\int_Q}f^2dxdt),
\end{array}
\eqno{(2.14)'}
$$
where $C$ depends only on $L_0,\ M, \ \Omega$.

\medskip
Now, we follow the standard method to provide a convergence proof
for the $\mathcal{K}$-approximate periodic solution $u^N(x,t)$ in
Proposition 2.3. Since some minor changes are needed, we give some
details. We first recall the energy estimate (Chapter 3 of [5]):
$$
\sup_{t\in[0,T]}\|u^N(\cdot,t)\|^2
+\displaystyle{\int^T_0}\|\nabla u^N(\cdot,t)\|^2 dt \leq
C(\|u^N(\cdot,0)\|^2+\int_Q f^2dxdt).
$$
By the estimate in (2.14)$'$, we get
$$
\begin{array}{ll}
\sup_{t\in[0,T]}\|u^N(\cdot,t)\|^2
+\displaystyle{\int^T_0}\|\nabla u^N(\cdot,t)\|^2 dt \leq C
(|a^N_I|^2 +\displaystyle{\int_Q}f^2dxdt),
\end{array}
\eqno{(2.15)}
$$
where $C$ depends only on $L_0, M, \Omega$, and where $u^N(x,t)$
is a solution of (2.7). We next fix
$a^N_I=a_I=(a_1,a_2,\cdots,a_{\mathcal{K}})$,  and also $\mathcal{K}$ that
satisfies the property in Remark 2.5. We compute
$$
\begin{array}{ll}
|u^N_j(t+\triangle t)-u^N_j(t)|&=|\displaystyle{\int^{t+\triangle
t}_t}(\sum B_{kj}(\tau)u^N_k (\tau)-f_j(\tau))d\tau|\\
&\leq |\displaystyle{\int^{t+\triangle
t}_t}(Lu^N,X_j)d\tau|+|\displaystyle{\int^{t+\triangle t}_t}
f_j(\tau)d\tau|\\
&\leq |\displaystyle{\int^{t+\triangle t}_t} \{ C_l(\varepsilon)
\|\nabla X_j\|^2+C_s(\varepsilon) \|\nabla X_j\|^2+
C_s(\varepsilon)\|\nabla u^N\|^2
+C\}d\tau| \\
&+(\displaystyle{\int^{t+\triangle t}_t}1\cdot
d\tau)^{\frac{1}{2}}(\displaystyle{\int^{t+\triangle t}_t}
f^2_j(\tau)d\tau
)^{\frac{1}{2}}\\
&\leq C_s(\varepsilon)+C_l(\varepsilon)\sqrt{\triangle t}.
\end{array}
$$
Now for any $\varepsilon'>0$, we can find $\varepsilon >0$ such
that $ C_s(\varepsilon)<\frac{\varepsilon'}{2}$. Then fixing such an
$\varepsilon$, we can find $\delta >0$ such that when $|\triangle
t|<\delta$, $C_l(\varepsilon)\sqrt{\triangle
t}<\frac{\varepsilon'}{2}$. Hence, when $|\triangle t|<\delta$,
$$
\begin{array}{ll}
|u^N_j(t+\triangle t)-u^N_j(t)|<\varepsilon'.
\end{array}
$$
Note that $\delta$ is independent of $N$. We proved that
$\{u^N_j(t)\}^{\infty}_{N=1}$ is equi-continuous. Since for any
$N>>1$, $\sum^{N}_{j=1}|u^N_j(t)|^2 \leq \|u^N(\cdot,t)\|^2 \leq
C$, $\{u^N_j(t)\}^{\infty}_{N=1}$ is uniformly bounded. Now by the
Ascoli-Arzela Theorem and the diagonal-element picking method, we
can find a subsequence $\{N_l\}$ such that for each $j$,
$u^{N_l}_j(t)\rightarrow u_j(t)$ uniformly over $[0,T]$.

\bigskip
{\bf Remark 2.6.}\ We note that for each $j$,
$\{u^N_j(t)\}^{\infty}_{N=1}$ is an equi-continuous family and
uniformly bounded. For each $j$ and for any $\varepsilon>0$, there
is a $\delta(j,\varepsilon)$ with $\delta$ depending only on
$j,L_0,M,\varepsilon$ such that when $|\triangle
t|<\delta(j,\varepsilon)$, we have $|u_j(t+\triangle
t)-u_j(t)|<\varepsilon$.\\

Next, for any fixed $m<N_l$, we have
$$
\begin{array}{ll}
\sum^m_{j=1}(u^{N_l}_j(t))^2\leq \sum
^{N_l}_{j=1}(u^{N_l}_j(t))^2=\|u^{N_l}(\cdot,t)\|^2\leq C.
\end{array}
$$
Letting $N_l\rightarrow \infty$, we get $\sum^m_{j=1}(u_j(t))^2
\leq C$, and thus $\sum^{\infty}_{j=1}(u_j(t))^2\leq C$. Let
$u(x,t)=\sum^{\infty}_{j=1}u_j(t)X_j(x)$. We have $u(x,t)\in
L^{\infty}(0,T;L^2(\Omega))$. Now,
 for $j\leq \mathcal{K}$,
$$
u_j(0)=\lim _{N_l\rightarrow \infty}u^{N_l}_j(0)=\lim
_{N_l\rightarrow \infty} a_j=a_j.
$$
For $ j> \mathcal{K}$,
$$
 u_j(0)=\lim _{N_l\rightarrow
\infty}u^{N_l}_j(0)=\lim _{N_l\rightarrow \infty}u^{N_l}_j(T)=u_j(T).
$$
Since $u^{N_l}(x,t),\nabla u^{N_l}(x,t)\in L^2(Q)$ with
$\|u^{N_l}(x,t)\|_{L^2(Q)}\leq C$, $\|\nabla
u^{N_l}(x,t)\|_{L^2(Q)}\leq C$, without loss of generality, we can
assume that $u^{N_l}(x,t) \rightarrow u^*(x,t)$ weakly in
$L^2(0,T;H^1_0(\Omega))$. (See page 54 of [5]).\\
\hspace*{0.5 cm} Let $\theta _j(t) \in C^{\infty} [0,T]$, and
$\Phi ^r(x,t)=
 \sum^r_{j=1} \theta _j(t)X_j(x)$, $N_l>r$. Since
$$
\begin{array}{ll}
 &\displaystyle{\int^T_0}(u^*(\cdot,t),\Phi ^r(\cdot,t))dt=
 \lim_{N_l\rightarrow \infty}
 \displaystyle{\int^T_0}(u^{N_l}(\cdot,t),\Phi ^r(\cdot,t))dt\\
 &= \lim_{N_l\rightarrow \infty}\displaystyle{\int^T_0}
 \sum^r_{j=1}u^{N_l}_j(t)\theta_j(t)dt=
 \displaystyle{\int^T_0}\sum^r_{j=1}u_j(t)\theta_j(t)dt\\
 &=\displaystyle{\int^T_0}(u(\cdot,t),\Phi ^r(\cdot,t))dt.
\end{array}
$$
Since the set of  all such $\Phi ^r(x,t)'s$ is dense in
$L^2(0,T;H^{1}_0(\Omega))$, we get $u(x,t)\equiv u^*(x,t)$ a.e. in
$L^2(0,T;H^{1}_0(\Omega))$. Hence, $u(x,t)\in L^\infty
(0,T;L^2(\Omega )) \cap L^2(0,T;H^1_0(\Omega ))$.

\medskip
Now we prove $u(x,t)$ is a weak solution.\\
\hspace*{0.5 cm} Let $\Phi ^r(x,t)=
 \sum^r_{j=1} \theta _j(t)X_j(x)$ be as above and $N_l>r$. We have
$$
 \lim_{N_l\rightarrow \infty}(u^{N_l}(\cdot,t),\Phi ^r(\cdot,t))=
 \lim_{N_l\rightarrow \infty}\sum^r_{j=1}u^{N_l}_j(t)\theta _j(t)
=\sum^r_{j=1}u_j(t)\theta _j(t)=(u(\cdot,t),\Phi ^r(\cdot,t)),
$$
$$
\lim_{N_l\rightarrow \infty}\displaystyle{\int^t_0}
 (u^{N_l}(\cdot,\tau),\Phi ^r_{\tau}(\cdot,\tau))d\tau =
 \displaystyle{\int^t_0} \lim_{N_l\rightarrow \infty}
 \sum^r_{j=1}u^{N_l}_j(\tau)\theta'_j(\tau)d\tau
 $$
 $$
=\displaystyle{\int^t_0}\sum^r_{j=1}u_j(\tau)\theta'_j(\tau)d\tau
=\displaystyle{\int^t_0}
 (u(\cdot,\tau),\Phi ^r_{\tau}(\cdot,\tau))d\tau,
$$
$$
\lim_{N_l\rightarrow
\infty}\displaystyle{\int^t_0}(Lu^{N_l}(\cdot,\tau),\Phi
^r(\cdot,\tau))d\tau=\displaystyle{\int^t_0}(Lu(\cdot,\tau),\Phi
^r(\cdot,\tau))d\tau.
$$
Therefore, we get
$$
(u(\cdot,t),\Phi ^r(\cdot,t))-\displaystyle{\int^t_0}
 (u(\cdot,\tau),\Phi ^r_{\tau}(\cdot,\tau))d\tau +
 \displaystyle{\int^t_0}(Lu(\cdot,\tau),\Phi^r(\cdot,\tau))d\tau
$$
$$
=(u(\cdot,0),\Phi^r(\cdot,0))+\displaystyle{\int^t_0}
(f(\cdot,\tau),\Phi^r(\cdot,\tau))d\tau .
$$
Since all such $\Phi ^r(x,t)'s$ are dense in $H^{1,1}_0(Q)$,  we
proved that  $u(x,t)$ is a weak solution of  (2.1). Note that
$u(x,t)\in\mathcal{S}_{\mathcal{K}}$. Summarizing the above, we
proved  the following:

\medskip
{\bf Theorem 2.7.}\ {\it There exists an integer ${\mathcal
K}_0\equiv{\mathcal K}_0(L_0,M,\Omega)$ $\geq 0$ such that for any
${\mathcal K}\geq{\mathcal K}_0$ and any initial value
$a_I=(a_1,a_2,\cdots,a_{\mathcal{K}})\in {\bf R}^{\mathcal{K}}$,
there is a unique solution to the following equation:
$$
\left\{
\begin{array}{ll}
\displaystyle{\frac{\partial u(x,t)}{\partial t}}+L_0u(x,t)+e(x,t)u(x,t)=f(x,t),
\;\;\;\;\;\;\;\;\;\;\;\;\;\;\;\;\;\;\;\;\;&
\mbox{in }\;\;Q,\\
 u(x,t)=0, & \mbox{on }\;\; \Sigma,\\
(u(\cdot,0),X_j)=a_j,&j\leq \mathcal{K},\\
(u(\cdot,0),X_j)=(u(\cdot,T),X_j),&j>\mathcal{K}.
\end{array}\right.
\eqno{(2.16)}
$$}

\medskip
{\bf Theorem 2.8.}\ {\it Let $u(x,t)$ be as in Theorem 2.7. Then
there is a constant $C$ depending only on $L_0,M,\Omega$ such that
$$
\begin{array}{ll}
\sup_{t\in[0,T]}\|u(\cdot,t)\|^2 +\displaystyle{\int^T_0}\|\nabla
u(\cdot,t)\|^2 dt \leq C (|a_I|^2 +\displaystyle{\int_Q}f^2dxdt).
\end{array}
\eqno{(2.17)}
$$}

\medskip
{\it Proof of Theorem 2.8.}\ By (2.15), we have
$$
\begin{array}{ll}
\sup_{t\in[0,T]}\|u^{N_l}(\cdot,t)\|^2
+\displaystyle{\int^T_0}\|\nabla u^{N_l}(\cdot,t)\|^2 dt \leq C
(|a_I|^2 +\displaystyle{\int_Q}f^2dxdt).
\end{array}
$$
Since $u^{N_l}(x,t) \rightarrow u^*(x,t)$ weakly in
$L^2(0,T;H^1_0(\Omega))$, we have
$$
\displaystyle{\int^T_0}\|\nabla u\|^2dt \leq \lim_{\overline{N_l
\rightarrow \infty}}\displaystyle{\int^T_0}\|\nabla u ^{N_l}\|^2dt
\leq C (|a_I|^2 +\displaystyle{\int_Q}f^2dxdt).\eqno{(2.18)}
$$
$$
\|u^{N_l}(\cdot,t)\|^2=\sum^{N_l}_{j=1}(u^{N_l}_j(t))^2 \leq C
(|a_I|^2 +\displaystyle{\int_Q}f^2dxdt).
$$
For any $m<N_l$, we have
$$
 \sum^m_{j=1} (u^{N_l}_j(t))^2\leq  C
(|a_I|^2 +\displaystyle{\int_Q}f^2dxdt).
$$
Letting $N_l\rightarrow\infty$, we get $\sum^m_{j=1}(u_j(t))^2\leq
C (|a_I|^2 +\displaystyle{\int_Q}f^2dxdt)$. Hence,
$$
\sum ^{\infty}_{j=1}(u_j(t))^2=\|u(\cdot,t)\|^2 \leq C (|a_I|^2
+\displaystyle{\int_Q}f^2dxdt).\eqno{(2.19)}
$$
By (2.18) and (2.19), we have (2.17). The proof is complete.
$\endpf$

\bigskip
{\it Proof of Theorem 2.7.}\ It suffices to prove the
uniqueness part of Theorem 2.7. Indeed, we need only to show that
the only solution $u(x,t)$ of (2.16) with $f=0,a_I=0$ is 0. For
this purpose, we first recall the energy estimate:
$$
\begin{array}{ll}
\sup_{t\in[0,T]}\|u(\cdot,t)\|^2 +\displaystyle{\int^T_0}\|\nabla
u(\cdot,t)\|^2 dt \leq C |a_{II}|^2,
\end{array}
\eqno{(2.20)}
$$
where
$$
\begin{array}{ll}
a_{II}=((u(\cdot,0),X_{{\mathcal{K}}+1}),(u(\cdot,0),X_
{{\mathcal{K}}+2}),\cdots\cdots).
\end{array}
$$
\hspace*{0.5 cm}Multiplying $X_j(x)$ to both sides of the first
equation of (2.16) and then integrating over $\Omega$, we get
$$
\begin{array}{ll}
\displaystyle{\frac{du_j(t)}{dt}} +(Lu,X_j)=0.
\end{array}
\eqno{(2.21)}
$$
\hspace*{0.5 cm}Write
$u^m(x,t)=\sum^m_{j={\mathcal{K}}+1}(u,X_j)X_j(x)=
\sum^m_{j={\mathcal{K}}+1}u_j(t)X_j(x),m>\mathcal{K}
+\mbox{1}$. Then $u^m(x,0)=u^m(x,T)$. \\
\hspace*{0.5 cm}Multiplying $2u_j(t)$ to (2.21) and summing up
with respect to $j$ from ${\mathcal{K}}+1$ to $m$, we get
$$
\begin{array}{ll}
\displaystyle{\frac{d(\|u^m(\cdot,t)\|^2)}{dt}} +2(L_0u,u^m)
+2(e(x,t)u,u^m)=0.
\end{array}
$$
By the same arguments as those in the proof of Lemma 2.1, we have
$$
\begin{array}{ll}
\displaystyle{\frac{(d\|u^m(\cdot,t)\|^2)}{dt}} +2\lambda
_{\mathcal{K}} \|u^m(\cdot,t)\|^2\leq C_s(\varepsilon)\|\nabla
u\|^2+C_l(\varepsilon)\|u\|^2.
\end{array}
$$
Using the Gronwall inequality and noticing that $u^m(x,0)=u^m(x,T)$, we get
$$
\begin{array}{ll}
(e^{2\lambda _{\mathcal{K}}T} -1)\|u^m(\cdot,0)\|^2
 &\leq C_s(\varepsilon)\displaystyle{\int^T_0}
 e^{2\lambda _{\mathcal{K}}t}
\|\nabla u\|^2dt+C_l(\varepsilon)\displaystyle{\int^T_0}
 e^{2\lambda_{\mathcal{K}}t}\|u\|^2dt \\
 &\leq C_s(\varepsilon)e^{2\lambda _{\mathcal{K}}T}|a_{II}|^2+
C_l(\varepsilon)\frac{(e^{2\lambda _{\mathcal{K}}T} -1)}
{2\lambda_{\mathcal{K}}}|a_{II}|^2.
\end{array}
$$
So
$$
\begin{array}{ll}
\|u^m(\cdot,0)\|^2\leq C_s(\varepsilon)\frac{e^{2\lambda
_{\mathcal{K}}T}}{e^{2\lambda _{\mathcal{K}}T} -1}|a_{II}|^2 +
C_l(\varepsilon)\frac{1}{2\lambda _{\mathcal{K}}}|a_{II}|^2.
\end{array}
$$
Letting $m\rightarrow\infty$, we obtain
$$
\lim_{m\rightarrow\infty}\|u^m(\cdot,0)\|^2=
\lim_{m\rightarrow\infty}\sum^m_{j={\mathcal{K}}+1}(u_j(0))^2
 =|a_{II}|^2 \leq C_s(\varepsilon)|a_{II}|^2 +
C_l(\varepsilon)\frac{1}{2\lambda _{\mathcal{K}}}|a_{II}|^2.
$$
We first choose $\varepsilon$ such that
$C_s(\varepsilon)<\frac{1}{4}$, then we choose ${\mathcal K}\ge
{\mathcal K}_{0}\gg \mbox{1}$ such that
$C_l(\varepsilon)\frac{1}{2\lambda _{\mathcal{K}}}<\frac{1}{4}$.
We have
$$
|a_{II}|^2\leq \frac{1}{4}|a_{II}|^2 + \frac{1}{4}|a_{II}|^2.
$$
It says $|a_{II}|^2\equiv 0$. By (2.20), we apparently get
$u(x,t)\equiv 0$. The proof is complete. $\endpf$

\medskip
{\it Proof of Theorem 1.2.} Theorem 1.2 follows directly from
Theorem 2.7-2.8. $\endpf$
\bigskip

\section{Existence of the solution to $(P)$}\ \hspace*{0.5 cm}
In this section, we give a proof of Theorem 1.3. Besides results
established in $\S 2$, another  main ingredient to be   used here
is the Carleman inequality for linear parabolic equation developed
in [2] [10] and [16], which in particular implies the  unique
continuation property for the solutions.




\medskip
First, by Theorem 2.7, there exists an integer ${\mathcal
K}_0\equiv{\mathcal K}_0$ $(L_0,M,\Omega)$ $\geq 0$ such that for
any ${\mathcal K}\geq{\mathcal K}_{0}$ and any initial value
$a_I=(a_1,a_2,\cdots,a_{\mathcal K})\in {\bf R}^{\mathcal{K}}$, we
have a unique solution $u(x,t)$ to the following equation:
$$
\left\{
\begin{array}{ll}
\displaystyle{\frac{\partial u(x,t)}{\partial t}}+L_0u(x,t)+e(x,t)u(x,t)=f(x,t),
\;\;\;\;\;\;&
\mbox{in }\;Q,\\
 u(x,t)=0, & \mbox{on }\;\; \Sigma,\\
(u(\cdot,0),X_j)=a_j,&j\leq \mathcal{K},\\
(u(\cdot,0),X_j)=(u(\cdot,T),X_j),&j>\mathcal{K}.
\end{array}\right.
\eqno{(3.1)}
$$

For the proof of Theorem 1.3, we need the following Lemma:

\medskip
{\bf Lemma 3.1.} {\it Let $e_m \in$ ${\mathcal M}_q$,
$a^m_I=(a^m_1,a^m_2,\cdots,a^m_{\mathcal K})\in{\bf
R}^{\mathcal{K}}$ with $m=1,2,\cdots $. And let $u_m\
(m=1,2,\cdots)$  be the solution of the following:
$$
\left\{
\begin{array}{ll}
\displaystyle{\frac{\partial u_m(x,t)}{\partial t}}+L_0u_m(x,t)+e_m(x,t)u_m(x,t)
=f_m(x,t), \;\;\;\;\;\;\;\;\;\;\;\;\;\;\;\;\;\;\;\;\;&
\mbox{in }\;\;Q,\\
 u_m(x,t)=0, & \mbox{on }\;\; \Sigma,\\
(u_m(\cdot,0),X_j)=a^m_j,&j\leq \mathcal{K},\\
(u_m(\cdot,0),X_j)=(u_m(\cdot,T),X_j),&j>\mathcal{K}.
\end{array}\right.
\eqno{(3.2)}
$$
Assume that $|a^m_I|\leq M_0$ with $M_0$ independent of the choice
of $m$. Suppose that $e_m\rightarrow e^*\in {\mathcal M}_q$ in the
weak-star topology of $L^{\infty}(0,T, L^q(\Omega))$, $a^m_j \rightarrow a^*_j$ for
$j=1,2,\cdots,\mathcal{K}$, and $f_m\rightarrow f^*$ in the
$L^2(Q)$-norm. Then there is a subsequence $\{m_k\}$ such that
$\{u_{m_k}\}$ converges in the weak $L^2$-topology to $u^*\in
C([0,T];L^2(\Omega))\cap L^2(0,T;H^1_0 (\Omega ))$ with
$$
\left\{
\begin{array}{ll}
\displaystyle{\frac{\partial u^*(x,t)}{\partial t}}+L_0u^*(x,t)
+e^*(x,t)u^*(x,t)=f^*(x,t),
\;\;\;\;\;\;\;\;\;\;\;\;\;\;\;\;\;\;\;\;\;&
\mbox{in }\;\;Q,\\
 u^*(x,t)=0, & \mbox{on }\;\; \Sigma,\\
(u^*(\cdot,0),X_j)=a^*_j,&j\leq \mathcal{K},\\
(u^*(\cdot,0),X_j)=(u^*(\cdot,T),X_j),&j>\mathcal{K}.
\end{array}\right.
\eqno{(3.3)}
$$
Moreover, write $u_{m_k}(x,t)=\sum^{\infty}_{j=1}
u^{m_k}_{j}(t)X_j(x)$ and
$u^*(x,t)=\sum^{\infty}_{j=1}u^*_{j}(t)X_j(x)$.  Then by choosing $m_k$ suitably, we have
$$
\lim_{k\rightarrow\infty}u^{m_k}_{j}(t)=u^*_j(t),\;\mbox{for any
}t\in [0,T], j,
$$
$ \mbox{and for any}\;\delta >0, u_{m_k}(x,t)\rightarrow
u^*(x,t)\mbox{ strongly in the }L^2([\delta,T]\times \Omega)$-norm. }

\medskip
{\it Proof of Lemma 3.1.} By the energy estimate in Theorem 2.7,
we have
$$
\begin{array}{ll}
\sup_{t\in[0,T]}\|u_m(\cdot,t)\|^2
+\displaystyle{\int^T_0}\|\nabla u_m(\cdot,t)\|^2 dt \leq C
(|a^m_I|^2 +\displaystyle{\int_Q}f^2_m dxdt)\leq C,
\end{array}
$$
where $a^m_I=(a^m_1,a^m_2,\cdots,a^m_{\mathcal{K}})$, $C$ depends
only on $L_0,M$, and $\Omega$. By the assumption in Lemma 3.1,
without loss of generality, we can assume $u_m(x,t)\rightarrow
u^*(x,t)$, $\nabla u_m(x,t)\rightarrow \nabla u^*(x,t)$ in the
weak  $L^2(Q)$-topology. Apparently, $u^*\in L^2(Q)\cap
L^2(0,T;H^1_0 (\Omega ))$.\\
\hspace*{0.5 cm} By Remark 2.6, $\{u^m_j(t)\}^{\infty}_{m=1}$ is
an equi-continuous family and uniformly bounded for each $j$. By
the diagonal-element picking method and the Ascoli-Arzela Theorem,
we can find a subsequence $\{m_k\}$ such that
$$
\begin{array}{ll}
\mbox{ for each }j,\;u^{m_k}_j(t)\rightarrow \widetilde{u^*_j}(t)
\mbox{ uniformly over }[0,T].
\end{array}
$$
Now, we let
$\widetilde{u^*}(x,t)=\sum^{\infty}_{j=1}\widetilde{u^*_j}(t)
X_j(x)$. Apparently, from the estimate:
$$
\begin{array}{ll}
\sum^N_{j=1}(u^{m_k}_j(t))^2 \leq C \mbox{ for any }N, \mbox{ it
follows that }\sum^{\infty}_{j=1}(\widetilde{u^*_j}(t))^2=
\|\widetilde{u^*}(\cdot,t)\|^2\leq C.
\end{array}
$$
Now for any $\theta_j(t)\in C^{\infty}[0,T]$, we have
$$
\begin{array}{ll}

\lim_{k\rightarrow\infty}\displaystyle{\int^T_0}
 (u_{m_k}(\cdot,t),X_j)\theta_j(t)dt =
\lim_{k\rightarrow\infty}\displaystyle{\int^T_0}
 u^{m_k}_j(t)\theta_j(t)dt \\
=\displaystyle{\int^T_0}
\widetilde{u^*_j}(t)\theta_j(t)dt=\displaystyle{\int^T_0}
(\widetilde{u^*}(\cdot,t),X_j)\theta_j(t)dt.
\end{array}
$$
On the other hand,
$$
\lim_{k\rightarrow\infty}\displaystyle{\int^T_0}(u_{m_k}
(\cdot,t),X_j)\theta_j(t)dt=\displaystyle{\int^T_0}(u^*
(\cdot,t),X_j)\theta_j(t)dt.
$$
Therefore,
$$
\begin{array}{ll}
\displaystyle{\int^T_0}(u^*
(\cdot,t)-\widetilde{u^*}(\cdot,t),X_j\theta_j(t))dt=0, & 0=
\displaystyle{\int^T_0}(u^* (\cdot,t)-\widetilde{u^*}(\cdot,t),
\sum^r_{j=1}X_j)\theta_j(t)dt\\   &=\displaystyle{\int_Q}(u^*
(x,t)-\widetilde{u^*}(x,t))( \sum^r_{j=1}X_j(x)\theta_j(t))dxdt\\
&=\displaystyle{\int_Q}(u^* (x,t)-\widetilde{u^*}(x,t)) \Phi
^r(x,t)dxdt.
\end{array}
$$
Since the set of functions with the form $\Phi ^r(x,t)=
 \sum^r_{j=1} \theta _j(t)X_j(x)$ is dense in
 $L^2(0,T;H^1_0(\Omega ))$, we
 conclude that $\widetilde{u^*}(x,t)=u^*(x,t)$ a.e. in $L^2(0,T;H^1_0(\Omega ))$.
 Therefore, $u^*(x,t)\in L^\infty (0,T;L^2(\Omega )) \cap L^2(0,T;H^1_0(\Omega ))$.\\
\hspace*{0.5 cm} Next, for the $\Phi ^r(x,t)$ defined above, by
the assumption,
$$
(u_m(\cdot,t),\Phi ^r(\cdot,t))-\displaystyle{\int^t_0}
 (u_m(\cdot,\tau),\Phi ^r_{\tau}(\cdot,\tau))d\tau +
 \displaystyle{\int^t_0}(L_0 u_m(\cdot,\tau),\Phi^r(\cdot,\tau))
 d\tau
$$
$$
+\displaystyle{\int^t_0}(e_{m}u_m(\cdot,\tau),\Phi^r(\cdot,\tau))
 d\tau
=(u_m(\cdot,0),\Phi^r(\cdot,0))+\displaystyle{\int^t_0}
(f_m(\cdot,\tau),\Phi^r(\cdot,\tau))d\tau .
$$
Now, as in Section 2, to show that $u^*(x,t)$ is the weak solution
of (3.3), it suffices to show that
$$
\begin{array}{ll}
\lim_{k\rightarrow\infty}\displaystyle{\int^T_0\int_{\Omega}}
e_{m_k}(x,\tau)u_{m_k}(x,\tau)
\Phi^r(x,\tau)dxd\tau=\displaystyle{\int^T_0\int_{\Omega}}
e^*(x,\tau)u^*(x,\tau)\Phi^r(x,\tau) dxd\tau
\end{array}
\eqno{(3.4)}
$$
for a certain subsequence  $\{m_k\}$. Notice that \\
for $j\leq \mathcal{K}$,
$$
\begin{array}{ll}
(u^*(\cdot,0),X_j)=u^*_j(0)=
\lim_{k\rightarrow\infty}u^{m_k}_j(0)=\lim_{k\rightarrow\infty}
a^{m_k}_j=a^*_j;
\end{array}
$$
for $j> \mathcal{K}$,
$$
\begin{array}{ll}
(u^*(\cdot,0),X_j)=
\lim_{k\rightarrow\infty}(u_{m_k}(\cdot,0),X_j)=
\lim_{k\rightarrow\infty}(u_{m_k}(\cdot,T),X_j)
=(u^*(\cdot,T),X_j).
\end{array}
$$
Hence, the proof of Lemma 3.1 will be complete if we can prove
(3.4).

Next, notice that
$$
\begin{array}{ll}
|\displaystyle{\int_Q} e_{m_k}u_{m_k} \Phi^rdxdt-
\displaystyle{\int_Q} e^*u^* \Phi^rdxdt| &\leq
|\displaystyle{\int_Q}(e_{m_k}- e^*)u^* \Phi^rdxdt| \\
&+|\displaystyle{\int_Q} e_{m_k}(u_{m_k}-u^*) \Phi^rdxdt|.
\end{array}
$$
Apparently, $|\displaystyle{\int_Q}(e_{m_k}- e^*)u^*
\Phi^rdxdt|\rightarrow 0$. Therefore, it suffices to prove the
following claim to complete the proof of Lemma 3.1.

\bigskip
{\bf Claim 3.2.} There is a subsequence $\{m_k\}$ such that
$|\displaystyle{\int_Q} e_{m_k}(u_{m_k}-u^*) \Phi^rdxdt|
\rightarrow 0$, and for any $\delta >0$,
$$
\begin{array}{ll}
u_{m_k}(x,t)\rightarrow u^*(x,t)\mbox{ strongly in the
}L^2([\delta,T]\times \Omega),\mbox{ as }k\rightarrow\infty.
\end{array}
$$

\medskip
{\it Proof of Claim 3.2.} Notice that
$$
 \left\{\begin{array}{ll}
\displaystyle{\frac{\partial(t^{\frac{2}{3}}u_m(x,t))}{\partial
t}}+L(t^{\frac{2}{3}} u_m(x,t))=t^{\frac{2}{3}}
f_m(x,t)+\frac{2}{3}t^{-\frac{1}{3}} u_m(x,t), \;\;\;\;\;&
\mbox{in }\;\;Q=\Omega \times (0,T),\\
 t^{\frac{2}{3}} u_m(x,t)=0, & \mbox{on }\;\; \Sigma=\partial
 \Omega \times (0,T),\\
 t^{\frac{2}{3}} u_m(x,t)|_{t=0}=0, & \mbox{in }\;\;\Omega.
\end{array}\right.
$$
By the high order energy estimates for parabolic equations
(page 59, Theorem 4.1 of [5]), we have
$$
\begin{array}{ll}
\sup_{t\in(0,T)} \|\nabla(t^{\frac{2}{3}} u_m(\cdot,t))\|^2+
 \displaystyle{\int_Q}|\partial_t( t^{\frac{2}{3}}
u_m(x,t))|^2dxdt &\leq  \displaystyle{\int_Q}|t^{\frac{2}{3}}
f_m(x,t)+\frac{2}{3}t^{-\frac{1}{3}} u_m(x,t)|^2dxdt\\
&\leq C.
\end{array}
$$
Then $\|t^{\frac{2}{3}} u_m(x,t)\|_{W^{1,1}(Q)}\leq C$ for all
$m$. By the Rellich lemma, there is a subsequence $\{m_k\}$ such
that $ t^{\frac{2}{3}} u_{m_k}(x,t)\rightarrow \overline{u}(x,t)$
strongly in $L^2(Q)$. Next,
$$
\begin{array}{ll}
&\displaystyle{\int_Q}\overline{u}(x,t)\Phi^r(x,t)
dxdt=\lim_{k\rightarrow \infty}
\displaystyle{\int_Q}t^{\frac{2}{3}} u_{m_k}(x,t)\Phi^r (x,t)dxdt\\
 &=\lim_{k\rightarrow \infty}\displaystyle{\int_Q}
 u_{m_k}(x,t)t^{\frac{2}{3}} \Phi^r(x,t)
dxdt=\displaystyle{\int_Q}t^{\frac{2}{3}}u^*(x,t)\Phi^r(x,t) dxdt.
\end{array}
$$
Since all such $\Phi^r(x,t)'s$ are dense in
$L^2(0,T;H^1_0(\Omega))$, $\overline{u}
(x,t)=t^{\frac{2}{3}}u^*(x,t)$ a.e. in $L^2(0,T;H^1_0(\Omega))$. Now,
$$
\begin{array}{ll}
&\displaystyle{\int^T_{\delta}\int_{\Omega}}|u_{m_k}-u^*|^2dxdt=
\displaystyle{\int^T_{\delta}\int_{\Omega}} t^{-4/3}
|t^{\frac{2}{3}}(u_{m_k}-u^*)|^2dxdt\\
&\leq C(\delta) \displaystyle{\int^T_{\delta}\int_{\Omega}}
|t^{\frac{2}{3}}(u_{m_k}-u^*)|^2dxdt\leq
C(\delta)\displaystyle{\int^T_{0}\int_{\Omega}}
|t^{\frac{2}{3}}(u_{m_k}-u^*)|^2dxdt  \rightarrow 0.
\end{array}
$$
We thus see, for any $ \delta >0$, that
$$
\begin{array}{ll}
u_{m_k}(x,t)\rightarrow u^*(x,t)\mbox{ strongly in the
}L^2([\delta,T]\times \Omega).
\end{array}
$$
Next,  by Claim 2.2,
$$
\begin{array}{ll}
|\displaystyle{\int _{\Omega}}e_{m_k}(u_{m_k}-u^*)\Phi^r dx| &\leq
C_s(\varepsilon)\displaystyle{\int
_{\Omega}}|e_{m_k}(u_{m_k}-u^*)^2|dx+C_l(\varepsilon)\displaystyle{\int
_{\Omega}}|e_{m_k}(\Phi^r)^2|dx \\
&\leq C_s(\varepsilon)(\displaystyle{\int
_{\Omega}}|\nabla(u_{m_k}-u^*)|^2dx +\displaystyle{\int
_{\Omega}}|(u_{m_k}-u^*)|^2dx)+C_l(\varepsilon) \\
&\leq C_s(\varepsilon)\displaystyle{\int
_{\Omega}}|\nabla(u_{m_k}-u^*)|^2dx + C_l(\varepsilon).
\end{array}
$$
Next,
$$
\begin{array}{ll}
|\displaystyle{\int^{\delta}_{0}\int_{\Omega}}e_{m_k}(u_{m_k}-u^*)
\Phi^r dxdt|&\leq
\displaystyle{\int^{\delta}_{0}}(C_s(\varepsilon)
\displaystyle{\int_{\Omega}}
|\nabla(u_{m_k}-u^*)|^2dx+C_l(\varepsilon))dt\\
&\leq C_s(\varepsilon)+\delta C_l(\varepsilon).
\end{array}
$$
Now, for any $\varepsilon'>0$, $\exists\varepsilon(\varepsilon')$
such that $C_s(\varepsilon)<\frac{\varepsilon'}{4}$. Then there is
a $\delta=\delta(\varepsilon,\varepsilon')$ such that $\delta
C_l(\varepsilon)<\frac{\varepsilon'}{4}$. We thus have
$$
\begin{array}{ll}
|\displaystyle{\int^{\delta}_{0}\int_{\Omega}}e_{m_k}(u_{m_k}-u^*)
\Phi^r dxdt|&\leq C_s(\varepsilon)+\delta C_l(\varepsilon)
 \leq \frac{\varepsilon'}{2}.
\end{array}
\eqno{(3.5)}
$$
Notice that
$$
\begin{array}{ll}
|\displaystyle{\int^{T}_{\delta}\int_{\Omega}}e_{m_k}(u_{m_k}-u^*)
\Phi^r dxdt|&\leq
C_l(\varepsilon'')\displaystyle{\int^{T}_{\delta}\int_{\Omega}}
|e_{m_k}(u_{m_k}-u^*)^2|dxdt\\
&+C_s(\varepsilon'')\displaystyle{\int^{T}_{\delta}\int_{\Omega}}
|e_{m_k}(\Phi^r)^2|dxdt\\
&\leq
C_l(\varepsilon'')\{C_l(\varepsilon''')\displaystyle{\int^{T}_{\delta}
\int_{\Omega}}|(u_{m_k}-u^*)^2|dxdt\\
&+ C_s(\varepsilon''')\displaystyle{\int^{T}_{\delta}
\int_{\Omega}}|\nabla(u_{m_k}-u^*)^2|dxdt\}+C_s(\varepsilon'').
\end{array}
$$
We have
$$
\begin{array}{ll}
\overline{\lim}_{k\rightarrow\infty}|\displaystyle{\int^{T}_{\delta}
\int_{\Omega}}e_{m_k}(u_{m_k}-u^*)\Phi^r dxdt| \leq
C_l(\varepsilon'')C_s(\varepsilon''')+C_s(\varepsilon'').
\end{array}
$$
For the $\varepsilon'$ as before , we can choose $\varepsilon''$
such that $C_s(\varepsilon'')<\frac{\varepsilon'}{4}$. Then for
this fixed $\varepsilon''$, there exists an $\varepsilon'''$ such that
$C_l(\varepsilon'')C_s(\varepsilon''')<\frac{\varepsilon'}{4}$.
Hence
$$
\begin{array}{ll}
\overline{\lim}_{k\rightarrow\infty}|\displaystyle{\int^{T}_{\delta}
\int_{\Omega}}e_{m_k}(u_{m_k}-u^*)\Phi^r dxdt| \leq
\frac{\varepsilon'}{2}.
\end{array}
\eqno{(3.6)}
$$
Since $\varepsilon'$ is arbitrary, by (3.5) and (3.6), we get
$|\displaystyle{\int_Q} e_{m_k}(u_{m_k}-u^*) \Phi^rdxdt|
\rightarrow 0$, as $k\rightarrow\infty$. The proof of Lemma 3.1 is
complete. $\endpf$

\bigskip {\it Proof of Theorem 1.3.} Let $d=\inf_{(e,a_I)\in
{\mathcal{M}}_q\times {\bf R}^{\mathcal{K}}}
\displaystyle{\int_{Q^{\omega}}}|u(e,a_I;x,t)-
\widetilde{u}|^{2}dxdt$. It is obvious that $d<\infty$. Thus there
exists a sequence $\{(e_m,a^m_I)\}_{m=1}^{\infty}$ such that
$$
d \leq \displaystyle{\int_{Q^{\omega}}}|u_m(e_m,a^m_I;x,t)-
\widetilde{u}|^{2}dxdt\leq d+\frac{1}{m}, \eqno{(3.7)}
$$
and
$$
\left\{
\begin{array}{ll}
\displaystyle{\frac{\partial u_m(e_m,a^m_I; x,t)}{\partial
t}}+L_0u_m(e_m,a^m_I;x,t)+e_m(x,t)u_m(e_m,a^m_I;x,t)=f(x,t), \;\;&
\mbox{in }\;Q,\\
 u_m(e_m,a^m_I;x,t)=0, & \mbox{on }\;\; \Sigma,\\
(u_m(e_m,a^m_I;\cdot,0),X_j)=a^m_j,&j\leq \mathcal{K},\\
(u_m(e_m,a^m_I;\cdot,0),X_j)=(u_m(e_m,a^m_I;\cdot,T),X_j),&j>\mathcal{K}.
\end{array}
\right. \eqno{(3.8)}
$$
In what follows, when there is no confusion of notation, we simply write
$u_m( x,t)$ for  $u_m(e_m,a^m_I; x,t)$.
By the definition of ${\mathcal M}_q$, there exists a subsequence
$\{m_k\}$ and $e^*\in{\mathcal M}_q$ such that
$$
e_{m_k}\rightarrow e^* \mbox{ in the weak star topology  as
}k\rightarrow\infty. \eqno{(3.9)}
$$
 By (3.7),
$$
\displaystyle{\int_{Q^{\omega}}}|u_m|^2dxdt\leq
\displaystyle{\int_{Q^{\omega}}}|\widetilde{u}|^2dxdt+C \leq C
\eqno{(3.10)}
$$
with $C$ independent of $m$.

\medskip
Next, we  prove the following claim:\\
\hspace*{0.5 cm} {\bf Claim 3.3.} {\it There is a constant $M_0$
such that $|a^m_I|^2\leq M_0$ for all $m$.}

\medskip
{\it Proof of Claim 3.3.} Suppose not . There is a subsequence
$\{m_k\}$ such that $\mu_{k}=|a^{m_k}_I|\rightarrow\infty$ as
$k\rightarrow\infty$. We write
$\hat{u}_{m_k}(x,t)=\frac{u_{m_k}(x,t)}{\mu_{k}}$. Then
$$
\left\{
\begin{array}{ll}
\displaystyle{\frac{\partial \hat{u}_{m_k}(x,t)}{\partial t}+L_0\hat{u}_{m_k}(x,t)
+e_{m_k}(x,t)\hat{u}_{m_k}(x,t) =\frac{f(x,t)}{\mu_{k}}},
\;\;\;\;\;\;\;\;\;\; &
\mbox{in }\;\;Q,\\
 \hat{u}_{m_k}(x,t)=0, & \mbox{on }\;\; \Sigma,\\
(\hat{u}_{mk}(\cdot,0),X_j)=\frac{a^{m_k}_j}{\mu_{k}},
&j\leq \mathcal{K},\\
(\hat{u}_{m_k}(\cdot,0),X_j)=(\hat{u}_{m_k}(\cdot,T),X_j),
&j>\mathcal{K}.
\end{array}
\right. \eqno{(3.11)}
$$
Apparently, $|\frac{a^{m_k}_I}{\mu_{k}}|^2=1$ for each $m_k$. After
passing to a subsequence, if necessary, we can assume that
$$
\begin{array}{ll}
 \frac{a^{m_k}_j}{\mu_{k}}\rightarrow \hat{a}_j \mbox{
with }
|\hat{a}_I|=1,\;\hat{a}_I=(\hat{a_1},\hat{a_2},\cdots,\hat{a}_{\mathcal{K}}),\;j=1,2,\cdots,
\mathcal{K}.
\end{array}
\eqno{(3.12)}
$$
Then by Lemma 3.1, there is a $\hat{u}\in C([0,T];L^2(\Omega))\cap
L^2(0,T;H^1_0 (\Omega ))$ such that
$\hat{u}_{m_k}\rightarrow\hat{u}$ in the weak $L^2$-topology,
$\mbox{and for any}\;\delta >0, \hat{u}_{m_k}(x,t)\rightarrow
\hat{u}(x,t)\mbox{ strongly in the }L^2([\delta,T]\times \Omega)$ topology.
 Also, all the other statements in Lemma 3.1 hold.\\
\hspace*{0.5 cm}By (3.7), we have
$$
\begin{array}{ll}
C&\geq
\displaystyle{\int_{Q^{\omega}}}|u_{m_k}(x,t)-\widetilde{u}|^2dxdt\\
&=\displaystyle{\int_{Q^{\omega}}}|\mu_{k}\hat{u}_{m_k}(x,t)-
\widetilde{u}|^2dxdt\\
&=(\mu_{k})^2\displaystyle{\int_{Q^{\omega}}}
|\hat{u}_{m_k}(x,t)- \frac{\widetilde{u}}{\mu_{k}}|^2dxdt\\
&\geq (\mu_{k})^2\displaystyle{\int_{[\delta,T]\times\omega}}
|\hat{u}_{m_k}(x,t)- \frac{\widetilde{u}}{\mu_{k}}|^2dxdt,\;\forall
\delta>0.
\end{array}
$$
Since $\mu_{k}\rightarrow\infty$, we see that
$$
\begin{array}{ll}
\lim_{k\rightarrow\infty}\displaystyle{\int_{[\delta,T]
\times\omega}}|\hat{u}_{mk}(x,t)-
\frac{\widetilde{u}}{\mu_{k}}|^2dxdt=0,\;\forall \delta>0.
\end{array}
$$
From the property that $\hat{u}_{m_k}(x,t)\rightarrow \hat{u}(x,t) \mbox{ strongly in
the }L^2([\delta,T]\times \Omega)$, it follows that
$$
\begin{array}{ll}
\displaystyle{\int_{[\delta,T]\times\omega}}
|\hat{u}(x,t)|^2dxdt=0,\;\forall \delta>0.
\end{array}
$$
By the unique  continuation property for solutions of the
parabolic equations, which is a  consequence the Carleman
inequality (see Page 430, Inequality (2) of [16]), we get
$\hat{u}(x,t)\equiv 0$ in $(0,T]\times \Omega$. Since  $\hat{u}\in
C([0,T],L^2(\Omega))$, we get $\hat{u}(x,0)\equiv 0$. On the other
hand, by (3.12), $(\hat{u}(\cdot,0),X_j)=\hat{a}_I$ with
$|\hat{a}_I|=1$. We see a contraction. The proof of Claim 3.3 is
complete. $\endpf$
\bigskip

Now, making use of  Claim 3.3, Lemma 3.1 and (3.9), we can assume that there
is a subsequence $\{m_k\}$ such that $a^{m_k}_I \rightarrow a^*_I$,
$\{u_{m_k}\}$ converges in the weak $L^2$-topology to $u^*\in
C([0,T];L^2(\Omega))\cap L^2(0,T;H^1_0 (\Omega ))$ with
$$
\left\{
\begin{array}{ll}
\displaystyle{\frac{\partial u^*(x,t)}{\partial t}}+L_0u^*(x,t)
+e^*(x,t)u^*(x,t)=f^*(x,t),
\;\;\;\;\;\;\;\;\;\;\;\;\;\;\;\;\;\;\;\;\;&
\mbox{in }\;\;Q,\\
 u^*(x,t)=0, & \mbox{on }\;\; \Sigma,\\
(u^*(\cdot,0),X_j)=a^*_j,&j\leq \mathcal{K},\\
(u^*(\cdot,0),X_j)=(u^*(\cdot,T),X_j),&j>\mathcal{K}.
\end{array}\right.
$$
By (3.7), we obtain
$$
\begin{array}{ll}
d\leq \displaystyle{\int_{Q^{\omega}}}|u^*(x,t)-
\widetilde{u}|^{2}dxdt\leq \lim_{\overline{k\rightarrow\infty}}
\displaystyle{\int_{Q^{\omega}}}|u_{m_k}(x,t)-
\widetilde{u}|^{2}dxdt\leq d.
\end{array}
$$
So
$$
\displaystyle{\int_{Q^{\omega}}}|u(e^*,a^*_I;x,t)-
\widetilde{u}|^{2}dxdt= \inf_{(e,a_I)\in {\mathcal{M}}_q\times
{\bf R}^{\mathcal{K}}}
\displaystyle{\int_{Q^{\omega}}}|u(e,a_I;x,t)-
\widetilde{u}|^{2}dxdt.
$$
The proof of Theorem 1.3 is complete. $\endpf$

\bigskip

{\it Proof of Corollary 1.4}: By Theorem 1.3, we can assume
that $k< {\mathcal K}_0$, where ${\mathcal K}_0$ is the same
integer  as in Theorem 1.2.  Notice that ${\mathcal S}_k
\subset {\mathcal S}_{{\mathcal K}_0}$. We keep the notation set up before.

Still  let $d=\inf_{(e,a_I)\in {\mathcal{M}}_q\times {\bf R}^{k},
u\in U(e,a_I;x,t)} \displaystyle{\int_{Q^{\omega}}}|u(x,t)-
\widetilde{u}|^{2}dxdt$ with $d<\infty$.  We then also have a
sequence of the pairs  $\{(e_m,a^m_I)\}_{m=1}^{\infty}$ and a
sequence of functions $u_m(e_m,a^m_I;x,t)$ such that
$$
d \leq \displaystyle{\int_{Q^{\omega}}}|u_m(e_m,a^m_I;x,t)-
\widetilde{u}|^{2}dxdt\leq d+\frac{1}{m}, \eqno{(3.13)}
$$
and
$$
\left\{
\begin{array}{ll}
\displaystyle{\frac{\partial u_m(e_m,a^m_I; x,t)}{\partial
t}}+L_0u_m(e_m,a^m_I;x,t)+e_m(x,t)u_m(e_m,a^m_I;x,t)=f(x,t), \;\;&
\mbox{in }\;Q,\\
 u_m(e_m,a^m_I;x,t)=0, & \mbox{on }\;\; \Sigma,\\
(u_m(e_m,a^m_I;\cdot,0),X_j)=a^m_j,&j\leq k,\\
(u_m(e_m,a^m_I;\cdot,0),X_j)=(u_m(e_m,a^m_I;\cdot,T),X_j),&j>k.
\end{array}
\right. \eqno{(3.14)}
$$
 Now, write $\widetilde{a}_I^{m}=(a_1^{m},\cdots,a_k^{m},
 \widetilde{a}_{k+1}^{m}, \widetilde{a}_{k+2}^{m},\cdots,\widetilde{a}_{{\mathcal K}_0}^{m})$
 with $\widetilde{a}_j^{m}=(u_m(e_m,a^m_I;x,0), X_j)=(u_m(e_m,a^m_I;x,T), X_j)$ for
$j=k+1,k+2,\cdots, {\mathcal K}_0$. Notice that
$$
\displaystyle{\int_{Q^{\omega}}}|u_m(e_m,a^m_I;x,t) |^2dxdt\leq
\displaystyle{\int_{Q^{\omega}}}|\widetilde{u}|^2dxdt+C \leq C.
\eqno{(3.15)}
$$
Since $u_m(e_m,a^m_I;x,t)\in {\mathcal S}_k \subset  {\mathcal
S}_{{\mathcal K}_0}$ for each $m$, making use of Lemma 3.1, we can
repeat the same argument as in the proof of Theorem 1.3 to show
that there is a subsequence $\{u_{m_l}(e_{m_l},a^{m_l}_I;x,t)\}$
of $\{u_m(e_m,a^m_I;x,t)\}$ such that
$\widetilde{a}^{m_l}_I\rightarrow a^{*}_I$ and $\{u_{m_l}\}$
converges in the weak $L^2$-topology to $u^*\in
C([0,T];L^2(\Omega))\cap L^2(0,T;H^1_0 (\Omega ))$ with
$$
\left\{
\begin{array}{ll}
\displaystyle{\frac{\partial u^*(x,t)}{\partial t}}+L_0u^*(x,t)
+e^*(x,t)u^*(x,t)=f(x,t),
\;\;\;\;\;\;\;\;\;\;\;\;\;\;\;\;\;\;\;\;\;&
\mbox{in }\;\;Q,\\
 u^*(x,t)=0, & \mbox{on }\;\; \Sigma,\\
(u^*(\cdot,0),X_j)=a^*_j,&j\leq {\mathcal K}_0,\\
(u^*(\cdot,0),X_j)=(u^*(\cdot,T),X_j),&j>{\mathcal K}_0.
\end{array}\right.
$$
Moreover, as in  Lemma 3.1, we also have

$$u^*_j(0)=
(u^*(\cdot,0),X_j)=
\lim_{m_l\rightarrow\infty}u^{m_l}_j(0)=
\lim_{m_l\rightarrow\infty}u^{m_l}_j(T)=u^*_j(T)
 \ \hbox {for }\ j>k.
$$
 Hence, $u^*\in {\mathcal S}_k$. Now, by the same argument as
 in the proof of Theorem 1.3, we conclude the proof of Corollary 1.4. $\endpf$
\bigskip

{\it Proof of Theorem 1.5}: Keep the notation as in Theorem 1.5.
Let $a_I^1, a_I^2\in {\mathbf R}^{\mathcal K}$ be such that
$$d=\displaystyle{\int_{Q^{\omega}}}|u(e,a^j_I;x,t)-
\widetilde{u}|^{2}dxdt= \inf_{a_I\in {\bf R}^{\mathcal{K}}}
\displaystyle{\int_{Q^{\omega}}}|u(e,a_I;x,t)-
\widetilde{u}|^{2}dxdt, \ \ j=1,\ 2.$$
For $\tau\in {\mathbf R}$, define
$$
I(\tau)=\displaystyle{\int_{Q^{\omega}}}|u(e,\tau a^1_I+(1-\tau)a^2_I;x,t)-
\widetilde{u}|^{2}dxdt.\ \ \hbox{Then}\
$$
$$
I(\tau)=\displaystyle{\int_{Q^{\omega}}}|\tau(u(e,a^1_I;x,t)-
\widetilde{u})+(1-\tau)(u(e,a^2_I;x,t)-
\widetilde{u})|^{2}dxdt.
$$
$$
=d(\tau^2+(1-\tau)^2)+2\tau(1-\tau)\displaystyle{\int_{Q^{\omega}}}
(u(e,a^1_I;x,t)-\widetilde{u})(u(e,a^2_I;x,t)-\widetilde{u})dxdt.
$$
Since $I(\tau)$ achieves its   minimum value at $\tau=0$, we have
$I'(0)=0$, from which the following follows:
$$
d=\displaystyle{\int_{Q^{\omega}}}(u(e,a^1_I;x,t)-\widetilde{u})(u(e,a^2_I;x,t)-\widetilde{u})dxdt.
$$
On the other hand, by the H\"older inequality, we have
$$
\left(\displaystyle{\int_{Q^{\omega}}}(u(e,a^1_I;x,t)-\widetilde{u})
(u(e,a^2_I;x,t)-\widetilde{u})dxdt\right)^2
$$
$$\le\displaystyle{\int_{Q^{\omega}}}|u(e,a^1_I;x,t)-
\widetilde{u}|^{2}dxdt\cdot \displaystyle{\int_{Q^{\omega}}}|u(e,a^2_I;x,t)-
\widetilde{u}|^{2}dxdt=d^2,
$$
with equality being held if and only if
$(u(e,a^1_I;x,t)-\widetilde{u})=C (u(e,a^2_I;x,t)-\widetilde{u})$
over $Q^{\omega}$ for a certain constant $C$. Apparently, this
implies that
$$
u(e,a^1_I;x,t)-\widetilde{u}= u(e,a^2_I;x,t)-\widetilde{u}\ \
\hbox{ over}\ Q^{\omega}.
$$
 By the Carleman inequality as mentioned in the proof of Theorem 1.3, we
 conclude that $u(e,a^1_I;x,t)= u(e,a^2_I;x,t)$ over $Q$. This then forces
 that they have the same initial value. In particular  we  conclude that
 $a_I^1=a^2_I$. The proof of Theorem 1.5 is complete.
$\endpf$
\bigskip

{\bf  Example 3.4}: Consider the following heat equation:

$$
y_t=\Delta y+cy+f(x,t), \ \ \ \; \; \; 0\le x\le 1, \;\; \ 0\le
t\le 1. \eqno (3.16)
$$
Here $c\in {\bf R}$ and $f\in L^2((0,1)\times (0,1))$. Notice that
$\{\frac{1}{\sqrt{2}}\sin(k\pi x)\}_{k=1}^{\infty}$ forms an
orthonormal basis of $L^2(0,1)$. Write
$f(x,t)=\sum_{k=1}^{\infty}f_k(t)\sin(k\pi x)$ and
$y=\sum_{k=1}^{\infty} a_k(t)\sin(k\pi x)$. Then we have
$a_k'(t)=-(k\pi)^2a_k(t)+ca_k(t)+f_k(t).$ Thus, we get
$$
a_k(1)\mbox{e}^{-c+(k\pi)^2}-a_k(0)= \int_0^1 f_k(t){\mbox
e}^{(-c+(k\pi)^2)t}dt.$$ Now, choose $c=(K\pi)^2$ and choose $f$
such that $\displaystyle{\int}_0^1 f_K(t)dt\not =0$. Then (3.16)
can never have a solution $y$, which is in the space ${\mathcal
S}_{K-1}$. However, in this case, for any
$b^K_I=(b_1,\cdots,b_{K})$, (3.16) does have a unique solution
$y(x,t)\in {\mathcal S}_{K}$ with $(y(\cdot,0), \frac{1}{\sqrt
2}\sin(j\pi x))=b_j$ for $j=1,2,\cdots, K$.

\par

\noindent {\bf Regular Mailing Address}: Room 5032, Dormitory 29,
Yuquan Campus, Zhejiang University, Hangzhou 310027, P. R.
China\\\\
\noindent {\bf E-mail Address}: leiling0810$@$yahoo.com.cn
\medskip

\end{document}